\documentclass[11pt]{article}

\usepackage{amsmath,amsthm}

\usepackage{amssymb,latexsym}

\usepackage{enumerate}

\newcommand{\dyw}{\mbox{\rm div}}

\newtheorem{tw}{Theorem}[subsection]
\newtheorem{lm}[tw]{Lemma}
\newtheorem{wn}[tw]{Corollary}
\newtheorem{stw}[tw]{Proposition}
\newenvironment{dow}{\it Proof.\rm}{\hfill $\Box$}

\theoremstyle{definition}
\newtheorem*{df}{Definition}
\newtheorem{uw}[tw]{Remark}
\newtheorem{prz}[tw]{Example}

\newcommand{\BN}{{\mathbb N}}
\newcommand{\BR}{{\mathbb R}}
\newcommand{\BX}{{\mathbb X}}

\newcommand{\FF}{{\mathcal{F}}}

\newcommand{\BB}{{\mathcal{B}}}

\newcommand{\MM}{{\mathcal{M}}}

\newcommand{\RR}{{\mathcal{R}}}
\newcommand{\EE}{{\mathcal{E}}}

\newcommand{\uAF}{\xrightarrow{uAF}}

\newcommand{\nsubsection}{\setcounter{equation}{0}\subsection}

\begin{document}
\title{Semilinear elliptic equations with Dirichlet operator
and singular nonlinearities}
\author{Tomasz Klimsiak }
\date{}
\maketitle

\begin{abstract}
In the paper we consider elliptic equations of the form
$-Au=u^{-\gamma}\cdot\mu$, where $A$ is the operator associated
with a regular symmetric Dirichlet form, $\mu$ is a positive
nontrivial measure and $\gamma>0$. We prove the existence and
uniqueness of solutions of such equations as well as some
regularity results. We also study  stability of solutions with
respect to the convergence of measures on the right-hand side of
the equation. For this purpose, we introduce some type of
functional convergence of smooth measures, which in fact is
equivalent to the quasi-uniform convergence of associated
potentials.
\end{abstract}
\noindent {\small\bf Mathematics Subject Classification (2010), 35J75, 60J45}.

\footnotetext{T. Klimsiak: Institute of Mathematics, Polish
Academy of Sciences, \'Sniadeckich 8, 00-956 Warszawa, Poland, and
Faculty of Mathematics and Computer Science, Nicolaus Copernicus
University, Chopina 12/18, 87-100 Toru\'n, Poland.

e-mail: tomas@mat.umk.pl; tel.: +48 566112951; fax: +48 56 6112987.}

\nsubsection{Introduction}

Let $E$ be a separable locally compact metric space,
$(\EE,D[\EE])$ be a regular symmetric Dirichlet form on
$L^2(E;m)$ and let $\mu$ be a nontrivial (i.e. $\mu(E)>0$)
positive Borel measure on
$E$. In the present paper we study elliptic equations of the form
\begin{equation}
\label{eq1.1}-Au=g(u)\cdot\mu, \quad u>0,
\end{equation}
where $A$ is the operator associated with $(\EE,D[\EE])$ and
$g:\BR^+\setminus\{0\}\rightarrow \BR^+$ is a continuous function satisfying
\begin{equation}
\label{eq1.2}
c_1\le g(u)\cdot u^\gamma\le c_2,\quad u>0
\end{equation}
for some $c_1,c_2, \gamma>0$. The model example of (\ref{eq1.1})
is the Dirichlet problem
\begin{equation}
\label{eqi.3} \left\{
\begin{array}{ll}-\Delta^{\alpha/2}u=u^{-\gamma}
\cdot\mu,\quad u>0,&\mbox{on } D,
\smallskip\\
 \,\,\, u=0&\mbox{on  } \BR^d\setminus D,
\end{array}
\right.
\end{equation}
where $\alpha\in (0,2]$, $\gamma >0$ and $D$ is a bounded  open
subset of $\BR^d$.

The paper consists of two parts. In the first part we address the
problem of existence, uniqueness and regularity of solutions of
(\ref{eq1.1}). In the second part we study stability of solutions
of (\ref{eq1.1}) with respect to the convergence of measures on
the right-hand side of the equation. The above problems were
treated in \cite{BrezisPonce} in case $A=\Delta$ and
\cite{BoccardoOrsina} in case $A$ is a uniformy elliptic
divergence form operator. Some different but related problems are
studied in \cite{MuratPoretta} in case $A$ is a Leray-Lions type
operator. The main aim of the present paper is to generalize the
results of \cite{BoccardoOrsina,BrezisPonce} to equations with
general (possibly nonlocal) operators corresponding to symmetric
Dirichlet forms. We also refine some results proved in
\cite{BoccardoOrsina,BrezisPonce,MuratPoretta} for equations with
local operators.

In the first part of the paper (Sections \ref{sec3} and
\ref{sec4}) we assume that $\mu$ belongs to the class $\RR$  of
smooth (with respect to capacity associated with $(\EE,D[\EE])$)
positive Borel measures on $E$ whose potential is $m$-a.e. finite
(see Section \ref{sec2} for details). It is known (see
\cite[Proposition 5.13]{KR:JFA}) that $\MM_{0,b}\subset\RR$, where
$\MM_{0,b}$ is the class of bounded smooth measures on $E$. In
general, the inclusion is strict. For instance, in case of
(\ref{eqi.3}), $\RR$ includes smooth Radon measures $\mu$ such
that $\int_D\delta^{\alpha/2}(x)\,\mu(dx)<\infty$, where
$\delta(x)= \mbox{dist}(x,\partial D)$ (see \cite[Example
5.2]{Kl:JFA1}).

The first difficulty we encounter when considering equation
(\ref{eq1.1}) is to define properly a solution. Here we give a
probabilistic definition of a solution of (\ref{eq1.1}) via the
Feynman-Kac formula. Namely, by a solution of (\ref{eq1.1}) we
mean a quasi-continuous function $u$ on $E$ such that $u>0$
quasi-everywhere (q.e. for short) with respect to the capacity Cap
naturally associated with  $(\EE,D[\EE])$ and for q.e. $x\in E$,
\[
u(x)=E_x\int_0^\zeta g(u)(X_t)\,dA^\mu_t.
\]
Here $\{(X_t)_{t\ge0},(P_x)_{x\in E}\}$ is a Hunt process with
life time $\zeta$  associated with the form $(\EE,D[\EE])$, $E_x$
is the expectation with respect to $P_x$ and $A^\mu$ is the
positive continuous additive functional in the Revuz
correspondence with $\mu$.

One reason for adopting here the probabilistic definition of a
solution is that unlike  problem (\ref{eqi.3}), for general $A$
one can not expect that $\inf_{x\in K} u(x)>0$ for every compact
$K\subset E$. Therefore the variational definition of a solution
considered in \cite{BoccardoOrsina} is not (at least directly)
applicable to general equations of the form (\ref{eq1.1}), because
we do not know whether $g(u)\cdot\mu$ is a Radon measure. The
probabilistic approach allows one to overcome the difficulty.
Another advantage lies in the fact that it allows one to cope with
the uniqueness problem.

In Section \ref{sec3} we prove several results on existence and
uniqueness of solutions of  (\ref{eq1.1}) and its generalization
(equation with mixed nonlinearities). It is worth pointing out
that the rather delicate problem of uniqueness (see \cite{Serrin})
was not addressed in \cite{BoccardoOrsina}.

Regularity of solutions of (\ref{eq1.1}) is studied in Section
\ref{sec4}. First, in Proposition \ref{stw1.2}, we generalize some
result proved in \cite{LiskevichSemenov}, and
then we use this generalization to prove that if $\mu$ is bounded
then for every $\gamma>0$ the function $u^{(\gamma+1)/{2}}$
belongs to the extended Dirichlet space $D_e[\EE]$ and there
exists $c(\gamma)>0$ such that
\[
\EE(u^{(\gamma+1)/2},u^{(\gamma+1)/2})\le c(\gamma)
c_2\|\mu\|_{TV},
\]
where $\|\mu\|_{TV}$ denotes the total variation norm of $\mu$. In
case of (\ref{eqi.3}) the above inequality gives the estimate of
$u^{(\gamma+1)/{2}}$ in the norm of the fractional Sobolev space
$H^{\alpha/2}_0(D)$.

In the second part of the paper (Sections \ref{sec5}--\ref{sec7}),
we study stability of solutions $u_n$ of the problems
\begin{equation}
\label{eqi.4} -Au_n=g(u_n)\cdot\mu_n, \quad u_n>0
\end{equation}
under different assumptions on the type of convergence of measures
$\mu_n$ and the limit measure $\mu$. We always assume that
$\{\mu_n\}$ is a sequence of smooth nontrivial Borel measures on
$E$ such that $\sup_{n\ge 1}\|\mu_n\|_{TV}<\infty$. As for $\mu$,
we distinguish two cases: $\mu\in\MM_{0,b}$, i.e. $\mu$ is bounded
and smooth, and $\mu\in\MM_b$, i.e. $\mu$ is a general bounded
Borel measure on $E$.

In Section \ref{sec5} we start with the study of the general case
$\mu\in\MM_b$. Our main result (Theorem \ref{tw1.5})  says that if
$\mu_n\rightarrow \mu$ vaguely then the sequence $\{\nu_n:=
g(u_n)\cdot\mu_n\}$ is tight in the vague topology and its every
limit point is a smooth measure. Moreover, if
$\nu_{n}\rightarrow\nu$ vaguely, then, up to a subsequence,
$u_{n}\rightarrow u$ $m$-a.e., where $-Au=\nu$.

In Section \ref{sec6} we address the case $\mu\in\MM_{0,b}$. We
first  introduce some type of convergence of smooth measures,
which is stronger then the vague and the narrow convergence. At
the same time, it is weaker then the convergence in the variation
norm, but nevertheless it preserves the smoothness property. This
new concept of convergence of $\{\mu_n\}$ to $\mu$ is defined via
some sort of uniform convergence of the sequence of additive
functionals $\{A^{\mu_n}\}$ to $A^{\mu}$, so we denote it by
$\xrightarrow{uAF}$. We prove (see Proposition \ref{stw1.3},
Proposition \ref{stw1.11}) that, up to a subsequence, the
convergence $\mu_n\xrightarrow{uAF} \mu$ is equivalent to the
quasi-uniform convergence of $\{u_n\}$ to $u$, where $u_n, u$ are
solutions of  the  problems
\begin{equation}
\label{eqi.5} -Au_n=\mu_n,\qquad -Au=\mu,
\end{equation}
respectively. Therefore it is possible to define the convergence
$\mu_n\xrightarrow{uAF} \mu$ analytically without recourse to the
notion of additive functional from the probabilistic potential
theory. Note that this analytical characterization of the
convergence $\mu_n\xrightarrow{uAF}\mu$ may be viewed as a
significant generalization of the stability result proved in
\cite{BrezisPonce}. Our main theorem on stability of (\ref{eqi.4})
(Theorem \ref{stw1.6}) says that if $\mu_n\xrightarrow{uAF}\mu$
then (up to a subsequence) $u_n\rightarrow u$ q.e., where $u$ is a
solution of (\ref{eq1.1}). We also show (see Proposition \ref{stw1.12}) that if
$\mu_n\xrightarrow{uAF}\mu$ then $\{\mu_n\}$ is locally
equidiffuse, which again confirms  the usefulness of our new
notion of the convergence of measures.

In Section  \ref{sec7} we return to the case of general measure
$\mu\in\MM_{b}$ but we assume that $E\subset\BR^d$ and $\mu$ is
approximated by mollification, i.e. $\mu_n=j_{1/n}*\mu$, where
$j_{1/n}$ is a mollifier. In our main result we also restrict our
attention to a class of operators including $\Delta^{\alpha/2}$,
$\alpha\in(0,2]$, as a special case. It is known that
$\mu\in\MM_b$ admits a unique decomposition
\[
\mu=\mu_c+\mu_d
\]
into the singular part $\mu_c$ with respect to Cap (the so-called
concentrated part) and an absolutely continuous part $\mu_d$ with
respect to  Cap (the so-called diffuse part). The case $\mu_c=0$
is covered by results of Section \ref{sec6}, because we show that
$j_{1/n}*\mu_d\xrightarrow{uAF} \mu_d$. The case $\mu_c\neq0$ is
much more involved, but can be handled by combining the results of
Section \ref{sec5} with those of Section \ref{sec6}. Before
describing our main result, we first make some comments on the
simplest case $A=\Delta$.

If $A=\Delta$ then from the inverse maximum principle (see
\cite{DupaignePonce}) one can deduce that the singular part
$\mu_c$ (with respect to the Newtonian capacity cap$_2$) is
responsible for explosions of the solution $u$ of (\ref{eq1.1}).
When $u$ explodes, $g(u)$ is formally equal to zero, so it seems
that in (\ref{eq1.1}) the absorption term $g$ forces some
reduction of $\mu_c$.  Several natural question arise here. The
first one is whether such reduction really occurs and whether the
whole singular part $\mu_c$ is reduced? Another question is
whether in investigating (\ref{eqi.3}) one should consider the
Newtonian capacity cap$_2$, or, maybe, it is better to consider
other capacities (for example $p$-capacities)? What happens if
$\Delta$ is replaced by a general Dirichlet operator $A$? In
\cite{BoccardoOrsina} partial answers to these questions are given
in case $A=\Delta$. Let $u_n$ be a solution of (\ref{eqi.4}) with
$A=\Delta$ and $\mu_n=g_n\cdot m$ with $\{g_n\}\subset
L^{\infty}(D;m)$, where $m$ is the Lebesgue measure on $D$. In
\cite{BoccardoOrsina} it is proved that if $\mu$ is orthogonal to
cap$_2$, (\ref{eq1.2}) is satisfied with $\gamma\ge1$ and
$\mu_n\rightarrow\mu$ in the narrow topology, then $u_n\rightarrow
0$. For $\gamma\in(0,1)$ similar result is proved in case $\mu$ is
orthogonal to the $p$-capacity with $p>2$ being the H\"older
conjugate to $q=\frac{d(\gamma+1)}{d-1+\gamma}$. Finally, let us
mention that the same problem of reduction of the singular part of
$\mu$ forced by absorption $g$ is considered in
\cite{MuratPoretta} in case $g$ is bounded and $A$ is a
Leray-Lions type operator (i.e. local operator).

In Theorem \ref{tw1.7}, the main result of Section \ref{sec7}, we
prove that in fact $g$ forces the reduction of the whole singular
part $\mu_c$ of $\mu$ for every $\gamma>0$. To be more specific,
we prove that if $u_{n}$ is a solution of (\ref{eqi.4}) with
$\mu_n=j_{1/n}\ast\mu$, then, up to a subsequence, $u_n\rightarrow
u$ $m$-a.e., where
\[
-Au=g(u)\cdot\mu_d,\quad u>0.
\]
The above result makes it legitimate to define solutions of
(\ref{eq1.1}) with bounded Borel measure $\mu$ as the solutions of
(\ref{eq1.1}) with $\mu$ replaced by $\mu_d$.
With this  definition, Theorem \ref{tw1.7} is the existence
theorem for (\ref{eq1.1}) with bounded Borel measure $\mu$.
Finally, note that $\mbox{Cap}=\mbox{cap}_2$ if $A=\Delta$ and
that the capacity cap$_2$ is absolutely continuous with respect to
the $p$-capacity for $p\ge 2$. Therefore in case $\gamma\in(0,1)$
our result strengthens the corresponding result from
\cite{BoccardoOrsina}. It should be stressed, however, that in
\cite{BoccardoOrsina} more general approximations $\{\mu_n\}$ of
$\mu$ are considered.

\nsubsection{Preliminaries} \label{sec2}

In the paper $E$  is a locally compact separable metric space and
$m$ is a positive Radon measure on $E$ such that Supp$[m]=E$. By
$(\EE,D[\EE])$ we denote a symmetric Dirichlet form on $L^2(E;m)$.
Recall that this means that
\begin{enumerate}
\item [($\EE.1$)] $\EE:D[\EE]\times D[\EE]\rightarrow\mathbb{R}$, where
$D[\EE]$ is a dense linear subspace of $L^2(E;m)$,
\item[($\EE.2$)] $\EE$ is bilinear,  $\EE(u,v)=\EE(v,u)$ and $\EE(u,u)\ge 0$,
$u,v\in D[\EE]$,
\item[($\EE.3$)] $\EE$ is closed, i.e. $D[\EE]$ equipped with the inner
product generated by the form $\EE_1$ is a Hilbert space (Here, as
usual, for $\alpha>0$ we set
$\EE_\alpha(u,v)=\EE(u,v)+\alpha(u,v)$, $u,v\in D[\EE]$, where
$(\cdot,\cdot)$ is the usual inner product in $L^2(E;m)$),
\item[($\EE.4$)] $\EE$ is Markovian, i.e. if $u\in D[\EE]$ then
$v:= (0\vee u)\wedge 1\in D[\EE]$
and $\EE(v,v)\le \EE(u,u)$.
\end{enumerate}

By Riesz's theorem, for every $\alpha>0$  and $f\in L^2(E;m)$
there exists a unique function $G_\alpha f\in L^2(E;m)$ such that
\[
\EE_\alpha(G_{\alpha}f,g)=(f, g),\quad g\in L^2(E;m).
\]
It is an elementary check that $\{G_\alpha,\, \alpha>0\}$ is a
strongly continuous contraction resolvent on $L^2(E;m)$. By
$\{T_t,\, t\ge 0\}$ we denote the associated semigroup and by
$(A,D(A))$ the operator generated by $\{T_t\}$. It is well known
(see \cite[Section 1.3]{Fukushima}) that $D(A)\subset D[\EE]$ and
\[
\EE(u,v)=(-Au,v),\quad u\in D(A),v\in D[\EE].
\]

In the whole paper we assume that $(\EE,D[\EE])$ is regular and  transient,
i.e.
\begin{enumerate}
\item [($\EE.5$)] (regularity) the space $D[\EE]\cap C_0(E)$ is dense
in $D[\EE]$ with respect to the $\EE_1$-norm and in $C_0(E)$ with
respect to the supremum norm,
\item [($\EE.6$)] (transience) there exists a strictly positive function
$g$ on $E$ such that
\[
\int_E |u(x)|g(x)\, m(dx)\le \|u\|_{\EE},\quad u\in D[\EE].
\]
where
\[
\|u\|_\EE=\sqrt{\EE(u,u)},\quad u\in D[\EE].
\]
\end{enumerate}

In the whole paper we fix $\varphi\in \mathcal{B}_b(E)$ such that
$\varphi>0$, $\int_E\varphi\,dm=1$, and we put $h=G_1 \varphi$,
$\pi=\varphi\cdot m$.

Given a Dirichlet form $(\EE,D[\EE])$ we define the capacity
Cap$:2^E\rightarrow \mathbb{R}^{+}$ as follows: for an open
$U\subset E$ we set
\[
\mbox{Cap}(U)=\EE_1(h_U,h_U),
\]
where $h_U$ is the reduced function of $h$ on $U$ (see
\cite[Chapter III]{MR}), and for arbitrary $A\subset E$ we set
\[
\mbox{Cap}(A)=\inf\{\mbox{Cap}(U); A\subset
U\subset E,\, U\,\,\, \mbox{open}\}.
\]

An increasing sequence $\{F_n\}$ of closed subsets of $E$ is
called nest if Cap$(E\setminus F_n)\rightarrow 0$ as $n\rightarrow
\infty$. A subset $N\subset E$ is called exceptional if
Cap$(N)=0$. We say that some property $P$ holds quasi everywhere
(q.e. for short) if a set for which it does not hold is
exceptional.

We say that a function $u$ defined q.e. on $E$ is quasi-continuous
if there exists a nest $\{F_n\}$ such that $u_{|F_n}$ is
continuous for every $n\ge 1$. It is known that each function
$u\in D[\EE]$ has a quasi-continuous $m$-version. From now on for
$u\in D[\EE]$ we always consider its quasi-continuous version.

A Borel measure $\mu$ on $E$ is called smooth if it does not
charge exceptional sets and there exists a nest $\{F_n\}$ such
that $|\mu|(F_n)<\infty,\, n\ge 1$. By $S$ we denote the set of
all positive smooth measures on $E$.

In the paper we also use the capacity CAP considered in
\cite[Chapter 2]{Fukushima}. We would like  to stress that the
notions of exceptional sets, quasi-continuity and smooth measures
defined with respect to  Cap and with respect to CAP are
equivalent. Therefore in the paper we may use the results of
\cite{Fukushima,MR} interchangeably.

By $S^{(0)}_0$ we denote the set of all measures $\mu\in S$ for
which there exists $c>0$ such that
\begin{equation}
\label{eq2.1} \int_E |u|\,d\mu\le c \sqrt{\EE(u,u)},\quad u\in
D[\EE],
\end{equation}

For a given Dirichlet form $(\EE,D[\EE])$ one can always define
the so-called extended Dirichlet space $D_e[\EE]$ as the set of
$m$-measurable functions on $E$ for which  there exists an
$\EE$-Cauchy sequence $\{u_n\}\subset D[\EE]$ convergent $m$-a.e.
to $u$ (the so-called approximating sequence). One can show that
for $u\in D_e[\EE]$ the limit $\EE(u,u)=\lim_{n\rightarrow
\infty}\EE(u_n,u_n)$ exists and does not depend on the
approximating sequence $\{u_n\}$ for $u$. Each element $u\in
D_e[\EE]$ has a quasi-continuous version. It is  known that
$(\EE,D[\EE])$ is transient iff the pair $(\EE,D_e[\EE])$ is a
Hilbert space. In the latter case for a given measure $\mu\in
S_{0}^{(0)}$ inequality (\ref{eq2.1}) holds for every $u\in
D_e[\EE]$.

In the sequel we say that $u:E\rightarrow\BR$ is measurable if it
is  universally measurable, i.e. measurable with respect to the
$\sigma$-algebra
\[
\mathcal{B}^{*}(E)=\bigcap_{\mu\in\mathcal{P}(E)}\mathcal{B}^\mu(E),
\]
where $\mathcal{P}(E)$ is the set of all probability measures on
$E$ and $\mathcal{B}^\mu(E)$ is the completion of $\mathcal{B}(E)$
with respect to the measure $\mu$.

By $\MM_b$ we denote the set of all bounded Borel measures on $E$
and by $\MM_{0,b}$ the subset of $\MM_b$ consisting of smooth
measures. We say that a positive Borel measure $\mu$ on $E$ is
nontrivial if $\mu(E)>0$.

Given a Borel measurable function $\eta$ on $E$ and a Borel
measure $\mu$ on $E$ we write
\[
(\eta,\mu)=\int_E\eta\,d\mu.
\]
By $u\cdot \mu$ w denote the Borel measure on $E$ defined as
\[
(f,u\cdot\mu)=(f\cdot u,\mu),\quad f\in \mathcal{B}(E)
\]
whenever the integrals exist.

Let us recall that for given measurable spaces $(S,\mathcal{S})$,
$ (T,\mathcal{T})$ a function $\kappa: S\times
\mathcal{T}\rightarrow \mathbb{R}^{+}\cup \{\infty\}$ is called a
kernel (from $S$ to $\mathcal{T}$) if for every $B\in \mathcal{T}$
the mapping $S\ni s\mapsto \kappa(s,B)$ is $\mathcal{S}$
measurable and for every fixed $s$ the mapping $\mathcal{T}\ni
B\mapsto \kappa(s,B)$ is a measure. Let us also recall that for
given measure $\mu$ on $\mathcal{S}$ and kernel $\kappa$ from $S$
to $\mathcal{T}$ one can consider its product $\mu\otimes\kappa$,
which by  definition  is a measure on
$\mathcal{S}\otimes\mathcal{T}$ defined as
\[
(\mu\otimes\kappa)(f)=\int_S\int_Tf(s,t)\,\kappa(s,dt)\,\mu(ds).
\]

With a regular symmetric Dirichlet form $(\EE,D[\EE])$ one can
associate uniquely  a Hunt process $\mathbb{X}=((X_t)_{,t\ge 0},
(P_x)_{x\in E}, (\FF_t)_{t\ge0},\zeta)$ (see \cite{Fukushima}). It
is related to $(\EE,D[\EE])$ by the formula
\[
T_tf(x)=E_xf(X_t),\quad t\ge 0,\quad m\mbox{-a.e.},
\]
where $E_x$ stands for the expectation with respect to the measure
$P_x$. For  $\alpha, t\ge 0$ and $f\in\mathcal{B}^{+}(E)$ we write
\[
R_\alpha f (x)=E_x\int_0^\zeta e^{-\alpha t}f(X_t)\,dt, \quad
p_tf(x)=E_xf(X_t),\quad x\in E.
\]
It is well known (see \cite[Section 5.1]{Fukushima} that for each
$\mu\in S$ there exists a unique positive continuous additive
functional $A^\mu$ in the Revuz duality with $\mu$. For $\mu\in S$
we write
\[
(R_\alpha \mu)(x)=E_x\int_0^\zeta e^{-\alpha t}\,dA^\mu_t,\quad
x\in E
\]
For simplicity we denote $R_0$ by $R$.

By $S_{00}^{(0)}$ we denote the set  of  all $\mu\in S_0^{(0)}$
such that $\mu(E)<\infty$ and $R\mu$ is bounded. We set
\[
\RR=\{\mu\in S:R\mu<\infty\mbox{ q.e.}\}.
\]
It is known (see \cite[Lemma 4.3, Proposition 5.13]{KR:JFA}) that
if $(\EE,D[\EE])$ is transient then $\MM_{0,b}^{+}\subset
\mathcal{R}$. Note that by \cite[Lemma 4.3]{KR:JFA}, if $\mu\in
\mathcal{R}$ then the function $R\mu$ is quasi-continuous. For an
equivalent definition of the class $\RR$ see remarks following
\cite[Lemma 3.1]{KR:CM}.

For a  Borel set $B$ we set
\[
\sigma_B=\inf\{t>0; X_t\in B\},\quad D_A=\inf\{t\ge 0; X_t\in B\},
\quad \tau_B=\sigma_{E\setminus B},
\]
i.e. $\sigma_B$ is the first hitting time of $B$, $D_A$ is the
first debut time of $B$ and $\tau_B$ is the first exit time of
$B$.

By $B^r$ we denote the set of regular points for $B$, i.e.
\[
B^r=\{x\in E; P_x(\sigma_B>0)=0\}.
\]

By $\mathcal{T}$ we denote the set of all  stopping times with
respect to the filtration $(\FF_t)_{t\ge 0}$ and by $\mathbf{D}$
the set of all measurable functions $u$ on $E$ for which the
family
\[
\{u(X_\tau),\, \tau\in\mathcal{T}\}
\]
is uniformly integrable with respect to the measure $P_x$ for q.e.
$x\in E$.

For a Borel measure $\mu$  on $E$ and $\alpha\ge 0$ by $\mu\circ
R_\alpha$ we denote the measure defined as
\[
(f,\mu\circ R_\alpha)=(R_\alpha f,\mu),\quad f\in\mathcal{B}(E)
\]
and by $P_\mu$ the measure
\[
P_\mu(A)=\int_E P_x(A)\, \mu(dx),\quad A\in \FF_\infty.
\]

Finally, let us recall that a positive measurable function $u$ on
$E$ is called excessive if
\[
p_tu\le u,\quad t\ge 0,
\]
and $u$ is called potential if it is  excessive and for every sequence
$\{T_n\}\subset\mathcal{T}$ such that $T_n\nearrow T\ge\zeta$,
\[
\lim_{n\rightarrow\infty}E_xu(X_{T_n})=0.
\]
for q.e. $x\in E$.

\nsubsection{Existence and uniqueness of solutions}
\label{sec3}

Let us recall that in the whole paper we assume that
$(\EE,D[\EE])$ satisfies $(\EE.1)$--$(\EE.6)$. As for $\mu$ and
$g$,  unless otherwise stated, in the paper we  assume that
$\mu\in S$ and $g:\BR^+\setminus\{0\}\rightarrow\BR^+$ is a
continuous function satisfying (\ref{eq1.2}). We also adopt the
convention that $g(0)=+\infty$, $g(+\infty)=0$.

\begin{uw}
\label{uw.ex} The class of forms satisfying $(\EE.1)$--$(\EE.6)$
is quite wide. For instance, it includes forms generated by
divergence form operators considered in \cite{BoccardoOrsina},
i.e. operators of the form
\[
Au(x)=\dyw(a(x)\nabla u(x)),\quad x\in D,
\]
where $D$ is a bounded open subset of $\BR^d$ and $a$ is a
symmetric bounded uniformly elliptic $d$-dimensional matrix. A
model example of nonlocal operator associated with form satisfying
$(\EE.1)$--$(\EE.6)$ is the fractional Laplacian
$\Delta^{\alpha/2}$ on $D$ with $\alpha\in(0,2)$. For the above
and some other interesting examples see, e.g., \cite[Chapter
1]{Fukushima}.
\end{uw}

\begin{df}
We say that a  measurable function $u:E\rightarrow\BR^+$ is a
solution of (\ref{eq1.1}) if
\begin{enumerate}
\item[(a)]$u$ is quasi-continuous and  $0<u(x)<\infty$ q.e.,
\item[(b)]for q.e. $x\in E$,
\begin{equation}
\label{eq3.1} u(x)=E_x\int_0^\zeta g(u(X_t))\,dA^{\mu}_t.
\end{equation}
\end{enumerate}
\end{df}

We will need  the following hypothesis:
\begin{enumerate}
\item[(H)] $g:\BR^+\setminus\{0\}\rightarrow\BR^+$ is nonincreasing.
\end{enumerate}

\subsubsection{Existence and uniqueness of solutions of
(\ref{eq1.1})}

We begin with a comparison and uniqueness result.

\begin{stw}
\label{stw1.1} Assume that $\mu_1, \mu_2$ are smooth measures such
that $0\le\mu_1\le\mu_2$ and $g_1, g_2
:\BR^+\setminus\{0\}\rightarrow\BR^+$ are measurable functions
such that $g_1(y)\le g_2(y)$ for $y>0$. Moreover, assume that
either $g_1$ or $g_2$ satisfies {\rm{(H)}}. If $u_1$ is a solution
of \mbox{\rm{(\ref{eq1.1})}} with data $g_1,\mu_1$ and $u_2$ is a
solution of \mbox{\rm{(\ref{eq1.1})}} with data
$g_2, \mu_2$ then $u_1\le u_2$ q.e.
\end{stw}
\begin{dow}
Without loss of generality we may assume that $g_2$ is
nonincreasing. By the Meyer-Tanaka formula, for q.e. $x\in E$ we
have
\begin{align*}
(u_1-u_2)^+(x)&\le
E_x\int_0^\zeta\mathbf{1}_{\{u_1-u_2>0\}}(X_t)
(g_1(u_1)(X_t)\,dA_t^{\mu_1}-g_2(u_2)(X_t)\,dA_t^{\mu_2})\\&
=E_x\int_0^\zeta\mathbf{1}_{\{u_1-u_2>0\}}(X_t)
g_1(u_1)(X_t)\,d(A_t^{\mu_1}-A_t^{\mu_2})\\&\quad
+E_x\int_0^\zeta\mathbf{1}_{\{u_1-u_2>0\}}(X_t)
(g_1(u_1)-g_2(u_1))(X_t)\,dA_t^{\mu_2}\\&\quad
+E_x\int_0^\zeta\mathbf{1}_{\{u_1-u_2>0\}}(X_t)
(g_2(u_1)-g_2(u_2))(X_t)\,dA_t^{\mu_2}.
\end{align*}
Since $\mu_1\le\mu_2$, $dA^{\mu_1}\le dA^{\mu_2}$ under
$P_x$ for q.e. $x\in E$ by the properties of the Revuz duality.
Therefore the first integral on the right-hand side of the above
equality is nonpositive. The second one is nonpositive since
$g_1\le g_2$ and $\mu_2\ge 0$. Finally, the third term is
nonpositive due to the fact that $g_2$ is nonincreasing and
$\mu_2\ge 0$. Hence $(u_1-u_2)^+(x)=0$ for q.e. $x\in E$,
which implies that $u_1\le u_2$ q.e.
\end{dow}

\begin{wn}
\label{wn1.1} Assume that $\mu\in S$ and $g$ satisfies
\mbox{\rm{(H)}}.  Then there exists  at most one solution of
{\mbox{\rm(\ref{eq1.1})}}.
\end{wn}

In what follows we will also need the following two hypotheses.
The first one was introduced by P.A. Meyer and is called Meyer's
hypothesis (L).
\begin{enumerate}
\item[(L)]For some (and hence for every) $\alpha>0,\,
\delta_{\{x\}}\circ R_\alpha \ll m$ for every $x\in E$, where
$\delta_{\{x\}}$ is the Dirac measure on $E$ concentrated at $x$.
\item [($\EE.7)$] For every nearly  Borel set $B$ such that
Cap$(B)>0$, $P_x(\sigma_B<\infty)>0$ for q.e. $x\in E$.
\end{enumerate}

\begin{uw}
(i) Hypothesis (L) is satisfied if there exists a Borel measurable
function $r_\alpha:E\times E\rightarrow \BR^{+}$ such that for
every $f\in L^2(E;m)$,
\[
R_{\alpha}f=\int_Ef(y)r_\alpha(\cdot,y)\, m(dy),\quad
m\mbox{-a.e.}
\]
It therefore clear that operators from Remark \ref{uw.ex} satisfy
(L).
\smallskip\\
(ii)  Hypothesis (L) is also called ``absolute continuity
condition". For equivalents for this property see \cite[Theorems
4.1.2, 4.2.4]{Fukushima}.
\end{uw}

\begin{uw}
\label{uw.exx}
\label{uw.str} Observe that if  ($\EE$.7) is satisfied then
$R\mu>0$ q.e. for every nontrivial $\mu\in S$. Indeed, let $F$ be
a quasi support of $A^\mu$. Then by \cite[Theorem
5.1.5]{Fukushima} it is also a quasi support of $\mu$. Since $\mu$
is nontrivial, Cap$(F)>0$. Therefore by ($\EE$.7),
$P_x(\sigma_F<\zeta)>0$ q.e. Since $F$ is  a quasi support of
$A^\mu$, $E_x\int_0^\zeta dA^\mu_t>0$ for q.e. $x\in F$. Hence for
q.e. $x\in E$ we have
\[
0<E_x E_{X_{\sigma_F}}\int_0^\zeta dA^\mu_t\le R\mu(x).
\]
\end{uw}

\begin{uw}
(i) It is known that ($\EE$.7) is satisfied if the form $(\EE,
D[\EE])$ is irreducible (see \cite[Theorem 4.7.1]{Fukushima}).
\smallskip\\
(ii) ($\EE$.7) is satisfied if the form $(\EE, D[\EE])$ satisfies
Meyer's hypothesis (L) and $r_\alpha(\cdot,\cdot)$ defined as
$r_\alpha(x,\cdot)\cdot m= \delta_{\{x\}}\circ R_\alpha$ is
strictly positive. Indeed, let $F$ be a closed set such that
Cap$(F)>0$. Then
\begin{equation}
\label{eq3.2} 0<\int_E r_\alpha(x,y)\,d\mu_F(y)=R_\alpha\mu_F(x) =
e^\alpha_F(x)=E_xe^{-\alpha\sigma_F},
\end{equation}
where $\mu_F$ is the smooth measure associated with the
equilibrium $e_F$ (see \cite[Theorem 2.1.5]{Fukushima}). The first
inequality in (\ref{eq3.2}) follows from the fact that $\mu_F$ is
nontrivial (since Cap$(F)>0$) and $r_\alpha(\cdot,\cdot)$ is
strictly positive. By (\ref{eq3.2}) we have
$P_x(\sigma_F<\infty)>0$ for q.e. $x\in E$.
\smallskip\\
(iii) From (ii) and Remark \ref{uw.exx} it follows that the
operators from Remark \ref{uw.ex} satisfy ($\EE$.7).
\end{uw}

\begin{stw}\label{stw2.b}
Assume that $\mu\in\mathcal{R}$ and
$g:\mathbb{R}\rightarrow\mathbb{R}^+$ is continuous and bounded.
Then if $g$ is nonincreasing or $(\EE,D[\EE])$ satisfies Meyer's
hypotheses \mbox{\rm{(L)}} then there exists a solution of the
equation
\begin{equation}
\label{eq2.mb}
-Au=g(u)\cdot\mu.
\end{equation}
Moreover, if $\mu$ is nontrivial, $g$ is strictly positive and
\mbox{\rm{($\EE$.7)}} is satisfied then $u>0$ q.e.
\end{stw}
\begin{dow}
First let us assume that $\mu\in S^{(0)}_{00}$.  Let us put
$V=(D_e[\EE],\|\cdot\|_\EE)$ and define $\Phi:V\rightarrow V$,
$\mathcal{A}: V\rightarrow V'$ by
\[
\Phi(u)=R(g(u)\cdot\mu),\quad\mathcal{A}u=-Au-g(u)\cdot\mu,\quad u\in V.
\]
That $\Phi(u)\in V$  follows from the fact that
$S^{(0)}_{00}\subset S^{(0)}_0$ and $R(S^{(0)}_0)\subset
D_e[\EE]$, while the fact that $\mathcal{A}u\in V'$ is a
consequence of the inclusion $S^{(0)}_0\subset V'$. Now we will
show some properties of the mappings $\mathcal{A}$, $\Phi$. If $g$
is nonincreasing then
\[
\langle \mathcal{A}u-\mathcal{A}v,u-v\rangle
=\|u-v\|_\EE-((g(u)-g(v))\cdot\mu,u-v)\ge \|u-v\|_\EE,\quad u,v\in V,
\]
where $\langle\cdot,\cdot\rangle$ is the duality pairing between
$V$ and $V'$. Thus $\mathcal{A}$ is strongly monotone, hence
coercive. It is also clear that $\mathcal{A}$ is hemicontinuous
and bounded. As for $\Phi$, let us first observe that $\|\Phi
(u)\|_\infty\le\|g\|_\infty\|R\mu\|_\infty$, $u\in V$. Moreover,
$\Phi$ is continuous. Indeed, let $u_n\rightarrow u$ and let
$v_n=\Phi(u_n)$, $v=\Phi(u)$. Then
\[
\|v-v_n\|_\EE=(v-v_n,(g(u)-g(u_n))\cdot\mu)
\le 2\|R\mu\|_\infty\|g\|_\infty\int_E|g(u)-g(u_n)|\,d\mu.
\]
Since $u_n\rightarrow u$ in $\EE$, there exists a subsequence
$(n')\subset (n)$ such that $u_{n'}\rightarrow u$ q.e. (see
\cite[Theorem 2.1.4]{Fukushima}). From this and the above
inequality it follows that $v_{n'}\rightarrow v$ in $\EE$. The
above argument shows that for every subsequence $(n')\subset (n)$
there exists a further subsequence $(n'')\subset (n')$ such that
$v_{n''}\rightarrow v$ in $\EE$, which implies that
$v_n\rightarrow v$ in $\EE$. Also observe that if $(\EE,D[\EE])$
satisfies Meyer's hypothesis (L), then $\Phi$ is compact. Indeed,
let $\{u_n\}\subset V$. Then
\begin{equation}
\label{eq2.b2} |v_n(x)-p_tv_n(x)|\le \|g\|_\infty
E_x\int_0^t\,dA^\mu_r,\quad t\ge 0
\end{equation}
for q.e. $x\in E$. By \cite[Theorem 2.2, Proposition 2.4]{Kl:MA1}
there exists a subsequence (still denoted by $(n)$) such that
$\{v_n\}$ is convergent q.e. Let $v=\lim_{n\rightarrow \infty}
v_n$. Then
\[
\|v-v_n\|_\EE=(v-v_n,(g(u)-g(u_n))\cdot\mu) \le
2\|g\|_\infty\int_E|v-v_n|\,d\mu,
\]
which converges to zero as $n\rightarrow\infty$. Now we may
conclude the existence result. In case  $g$ is nonincreasing the
existence of a solution of (\ref{eq2.mb}) follows from
\cite[Corollary II.2.2]{Showalter}. If $(\EE,D[\EE])$ satisfies
Meyer's hypotheses (L) then the existence follows by the Schauder
fixed point theorem.

Now we turn to the the general case where $\mu\in\mathcal{R}$.
There exists a nest $\{F_n\}$ such that
$\mathbf{1}_{F_n}\cdot\mu\in S^{(0)}_{00}$, $n\ge1$ (see
\cite[Section 2.2]{Fukushima}). By what has already been proved,
for each $n\ge1$ there exists a solution $u_n\in V$ of the
equation
\[
-Au_n=g(u_n)\cdot\mu_n.
\]
By the definition of a solution,
\[
u_n(x)=E_x\int_0^\zeta g(u_n)\mathbf{1}_{F_n}(X_t)\,dA^\mu_t
\]
for q.e. $x\in E$. Since $\{F_n\}$ is a nest,
$\mathbf{1}_{F_n}(X_t)\rightarrow 0$, $t\in [0,\zeta)$, $P_x$-a.s.
for  q.e. $x\in E$ (see \cite[Proposition IV.5.30]{MR}). If $g$ is
nonincreasing then by Proposition \ref{stw1.1} the sequence
$\{u_n\}$ is nondecreasing. Therefore $u:=\lim_{n\rightarrow
\infty} u_n$ is a solution of (\ref{eq2.mb}). If $(\EE,D[\EE])$
satisfies Meyer's hypotheses (L) then by (\ref{eq2.b2}), which
holds with $v_n$ replaced by $u_n$, and by \cite[Theorem 2.2,
Propositions 2.4 and 4.3]{Kl:MA1}, there exists a subsequence
$(n')\subset (n)$ such that $\{u_{n'}\}$ is convergent q.e. It is
clear that $u:=\lim_{n'\rightarrow \infty}u_{n'}$ is a solution of
(\ref{eq2.mb}). The second assertion of the theorem follows
immediately from the assumptions and Remark \ref{uw.str}.
\end{dow}

\begin{lm}
\label{lm3.cap} Let $\mu\in \RR$ and let $u$ be  defined as
\[
u(x)=E_x\int_0^\zeta dA^\mu_t,\quad x\in E.
\]
Then
\[
\lim_{n\rightarrow\infty}\mbox{\rm CAP}(\{u>n\})\rightarrow 0.
\]
\end{lm}
\begin{dow}
Let $A_n=\{u>n\}$. If $\sigma_{A_n}<\infty$ then
$\sigma_{A_n}<\zeta$.  Therefore by the Markov property and the
fact that $\mu\in\RR$, for q.e. $x\in E$ we have
\[
P_x(\sigma_{A_n}<\infty)\le P_x(u(X_{\sigma_{A_n}\wedge\zeta})\ge
n) \le n^{-1}E_x\int_0^\zeta dA^\mu_t,
\]
which converges to zero as $n\rightarrow\infty$. Therefore
applying \cite[Corollary 4.3.1]{Fukushima} we get the desired
result.
\end{dow}

\begin{tw}
\label{tw1.1} Assume that $(\EE,D[\EE])$  satisfies $(\EE.7)$,
$\mu\in\mathcal{R}$ is nontrivial and $g$ satisfies
$\mbox{\rm{(H)}}$. Then there exists a solution of
\mbox{\rm(\ref{eq1.1})}.
\end{tw}
\begin{dow}
By Corollary \ref{wn1.1} and Proposition \ref{stw2.b}, for every
$n\ge1$ there exists a unique solution $u_n$ of the problem
\begin{equation}
\label{eq1.3} -Au_n=g_n(u_n)\cdot\mu,\quad u_n>0
\end{equation}
with $g_n(u)=g(u+\frac1n)$, $u>0$ and $g_n(u)=g(\frac1n)$, $u\le
0$. By Proposition \ref{stw1.1}, $\{u_n\}$ is nondecreasing. Hence
$u_1\le u_n$ for $n\ge 1$. Since $(\EE,D[\EE])$ satisfies
($\EE.7$) and $\mu$ is nontrivial, $u_1>0$ q.e. Hence $u_n>u_1>0$,
$n\ge 1$ q.e. Put $u=\limsup_{n\rightarrow\infty}u_n>0$. Then
$u>0$ q.e. By the definition of a solution of (\ref{eq1.3}),
\begin{equation}
\label{eq1.4} u_n(x)=E_x\int_0^\zeta g_n(u_n(X_t))\,dA_t^\mu
\end{equation}
for q.e. $x\in E$. By the Meyer-Tanaka formula and (\ref{eq1.2}),
\[
u_n^{\gamma+1}(x)\le(\gamma+1)E_x\int_0^\zeta g_n(u_n)u_n^\gamma
(X_t)\,dA_t^{\mu_\cdot}\le(\gamma+1)c_2 E_x\int_0^\zeta\,dA_t^\mu.
\]
Hence
\[
u_n^{\gamma+1}(x)\le(\gamma+1)c_2 E_x\int_0^\zeta\,dA_t^\mu<\infty
\]
for q.e. $x\in E$. From the above inequality we conclude that $u$
is a potential and $u\in \mathbf{D}$. Let $\tau_k=\tau_{G_k}$,
$G_k=\{u_1\ge k^{-1}\}$. Observe that for every $x\in G_k$,
\[
g(u_n(x)+\frac 1 n)\le g(u_1(x)+\frac 1 n)
\le\frac{c_2}{u_1^\gamma (x)}\le c_2 k^\gamma.
\]
Therefore by the Lebesgue dominated convergence theorem,
\[
E_x\int_0^{\tau_k}g_n(u_n)(X_t)\,dA_t^\mu\rightarrow
E_x\int_0^{\tau_k}g(u)(X_t)\,dA_t^\mu
\]
as $n\rightarrow\infty$. Since for each $k\ge 1$,
\[
u_n(x)=E_xu_n(X_{\tau_k})+E_x\int_0^{\tau_k}g_n(u_n)(X_t)\,dA^{\mu}_t
\]
for q.e. $x\in E$, it follows that
\[
u(x)=E_xu(X_{\tau_k})+E_x\int_0^{\tau_k}g(u)(X_t)\,dA_t^\mu
\]
for q.e. $x\in E$. Since $u$ is a potential, from Lemma
\ref{lm3.cap} and \cite[Lemma 5.1.6]{Fukushima} it follows that
$\lim_{k\rightarrow \infty}\tau_k\ge \zeta$. Therefore letting
$k\rightarrow\infty$ in the above equation we conclude that
(\ref{eq3.1}) is satisfied for q.e. $x\in E$.
\end{dow}

\begin{tw}
\label{tw1.2} Assume that $(\EE,D[\EE])$  satisfies $(\EE.7)$ and
Meyer's hypothesis \mbox{\rm(L)} and that $\mu\in\mathcal{R}$ is
nontrivial. Then there exists a solution of
$\mbox{\rm{(\ref{eq1.1})}}$.
\end{tw}
\begin{dow}
By Proposition \ref{stw2.b}, for every $n\ge1$ there exists a
solution $u_n$ of (\ref{eq1.3}) with $g_n(u)=g(u+\frac1n)$ for
$u>0$. By (\ref{eq1.2}) and Proposition \ref{stw1.1},
\[
v_n\le u_n\le w_n,\quad n\ge1\quad\mbox{q.e.},
\]
where $v_n,w_n$ are solutions of the problems
\begin{equation}
\label{eq3.7} -Av_n=c_1(v_n+\frac 1 n)^{-\gamma}\cdot\mu,\quad
v_n>0,\qquad -Aw_n=c_2(w_n+\frac1 n)^{-\gamma}\cdot\mu, \quad
w_n>0.
\end{equation}
Hence
\[
g(u_n+\frac 1 n)\le c_2(u_n+\frac 1 n)^{-\gamma}\le c_2(v_n+\frac
1 n)^{-\gamma}\quad\mbox{q.e.}
\]
Let $v,w$ be solutions of the problems
\[
-Av=c_1v^{-\gamma}\cdot\mu,\quad v>0,\qquad
-Aw=c_2w^{-\gamma}\cdot\mu,\quad w>0.
\]
From the proof of Theorem \ref{tw1.1} it follows that $\{v_n\}$
converges q.e. to $v$. Hence
\[
c_2(v_n+\frac 1 n)^{-\gamma}(X)\rightarrow c_2v^{-\gamma}(X),
\quad P_x\otimes dA^\mu\mbox{-a.s.}
\]
for q.e. $x\in E$, where $P_x\otimes dA^\mu$ is the product of the
measure  $P_x$ and the kernel $dA^\mu$ from $\Omega$ to
$\mathcal{B}(\mathbb{R}^+)$. Moreover,
\[
v_n(x)=E_x\int_0^\zeta c_2(v_n+\frac1n)^{-\gamma}(X_t)\,dA_t^\mu
\rightarrow E_x\int_0^\zeta c_2v^{-\gamma}(X_t)\,dA_t^\mu=v(x)
\]
for q.e. $x\in E$, which implies that the family
$\{c_2(v_n(X)+\frac1 n)^{-\gamma}\}$ is uniformly integrable with
respect to the measure $P_x\otimes dA^\mu$ for q.e. $x\in E$. From
this we conclude that
\[
\lim_{t\rightarrow 0^+}
\sup_{n\ge1}E_x\int_0^tg_n(u_n)(X_r)\,dA_r^\mu=0
\]
for q.e. $x\in E$. Therefore
\begin{equation}
\label{eq1.5a} \lim_{t\rightarrow 0^+}\sup_{n\ge
1}|u_n(x)-p_tu_n(x)|= \lim_{t\rightarrow 0^+}\sup_{n\ge
1}E_x\int_0^tg_n(u_n)(X_r)\,dA_r^\mu=0
\end{equation}
for q.e. $x\in E$. Since $u_n\le w$ for $n\ge1$, it follows from
\cite[Theorem 2.2, Propositions 2.4 and 4.3]{Kl:MA1} that there
exists a subsequence (still denoted by $(n)$) such that $\{u_n\}$
converges q.e. The rest of the proof runs as the proof of Theorem
\ref{tw1.1}.
\end{dow}

\subsubsection{Existence and uniqueness of solutions with mixed
nonlinearities}

In this subsection we study problems of the form
\begin{equation}
\label{eq1.6} -Au=(g(u)+h(u))\cdot\mu,\quad u>0.
\end{equation}

\begin{tw}
\label{tw1.3} Assume that $(\EE,D[\EE])$ satisfies $(\EE.7)$,
$\mu\in\RR$ is nontrivial,  $g,\, h$  satisfy \mbox{\rm(H)} and
$h:\BR^+\setminus \{0\}\rightarrow\BR^+$ is a continuous
 function  such that
\begin{equation}
\label{eq3.8} c_1\le h(s)\cdot s^\beta\le c_2,\quad s>0
\end{equation}
for some $\beta>0$.  Then there exists a unique solution $u$ of
problem \mbox{\rm(\ref{eq1.6})}. Moreover,
\begin{equation}
\label{eq1.7}
u\le\frac{c_2}{c_1}(2^\gamma v+2^\beta w),
\end{equation}
where $v,w$ are solutions of the problems
\[
-Av=c_1v^{-\gamma}\cdot\mu,\quad v>0,\qquad
-Aw=c_1w^{-\beta}\cdot\mu,\quad w>0.
\]
\end{tw}
\begin{dow}
Uniqueness follows from Proposition \ref{stw1.1}. To prove the
existence of solutions, let $u_n$ denote  the solution of the
problem
\begin{equation}
\label{eq.1001}
-Au_n=(g_n(u_n)+h_n(u_n))\cdot\mu,\quad u_n>0
\end{equation}
with $g_n(u)= g(u+\frac1n)$, $h_n(u)=h(u+\frac1n)$ for $u>0$. By
Proposition \ref{stw1.1}, $\{u_n\}$ is nondecreasing and
\begin{equation}
\label{eq1.8} v_n\le u_n,\quad w_n\le u_n\quad\mbox{q.e.},
\end{equation}
where $v_n,w_n$ are solutions of (\ref{eq3.7}). Therefore for each $n\ge1$,
\[
v_n+w_n\le 2 u_n\quad\mbox{q.e.}
\]
By Proposition \ref{stw1.1} the sequences $\{w_n\},\{v_n\}$ are
also nondecreasing.  Furthermore,
\begin{align}
\label{eq1.9} \nonumber g(u_n+\frac1n)+h(u_n+\frac 1 n)&\le
c_2(u_n+\frac 1 n)^{-\gamma} + c_2(u_n+\frac1n)^{-\beta}\\
&\le c_2(\frac12 w_n+\frac12v_n+\frac1n)^{-\gamma}
+c_2(\frac12w_n+\frac12v_n+\frac1n)^{-\beta}\nonumber\\
&\le c_2 2^\gamma(v_n+\frac1n)^{-\gamma}+
c_22^\beta(w_n+\frac1n)^{-\beta}.
\end{align}
From the proof of Theorem \ref{tw1.1} it follows that the
sequences $\{(v_n+\frac1n)^{-\gamma}(X)\}$ and
$\{(w_n+\frac1n)^{-\gamma}(X)\}$ are uniformly integrable with
respect to the measure $P_x\otimes dA^\mu$. Let
$u=\limsup_{n\rightarrow\infty}u_n$. By the definition of a
solution of (\ref{eq.1001}),
\begin{equation}
\label{eq1.10}
u_n(x)=E_x\int_0^\zeta(g_n(u_n)(X_t)+h_n(u_n)(X_t))\,dA_t^\mu
\end{equation}
for q.e. $x\in E$. By (\ref{eq1.9}) the sequence
$\{(g_n(u_n)(X)+h_n(u_n)(X)\}$ is uniformly integrable with
respect to the measure $P_x\otimes dA^\mu$. Therefore letting
$n\rightarrow\infty$ in (\ref{eq1.10}) we get
\[
u(x)=E_x\int_0^\zeta(g(u)(X_t)+h(u)(X_t))\,dA_t^\mu.
\]
Inequality (\ref{eq1.7}) follows easily from (\ref{eq1.9}.)
\end{dow}

\begin{tw}
\label{tw1.4} Assume that $(\EE,D[\EE])$ satisfies $(\EE.7)$ and
Meyer's  hypothesis {\mbox{\rm(L)}}, $\mu\in\RR$ is nontrivial and
$h:\BR^+\setminus\{0\}\rightarrow\BR^+$ is a continuous function
satisfying \mbox{\rm(\ref{eq3.8})} for some $\beta>0$. Then there
exists a solution of \mbox{\rm(\ref{eq1.6})} such that estimate
\mbox{\rm(\ref{eq1.7})} holds true.
\end{tw}
\begin{dow}
In the proof of Theorem \ref{tw1.3} monotonicity of $g,h$ was used
only to prove q.e. convergence of $\{u_n\}$. As in the proof of
Theorem \ref{tw1.3} we show that the sequence
$\{(g_n(u_n)(X)+h_n(u_n)(X)\}$ is uniformly integrable with
respect to the measure $P_x\otimes dA^\mu$. Therefore
(\ref{eq1.5a}) is satisfied, which when combined with
\cite[Theorem 2.2, Propositions 2.4 and 4.3]{Kl:MA1} implies that
$\{u_n\}$ has a subsequence convergent q.e.
\end{dow}

\nsubsection{Regularity of solutions} \label{sec4}

\begin{df}
\label{df1.1} We say that a sequence $\{u_n\}$ of  measurable
functions is convergent quasi-uniformly to a function $u$ if for
every $\varepsilon>0$,
\begin{equation}
\label{eq1.14}
\lim_{n\rightarrow\infty}\mbox{CAP}(\{|u_n-u|>\varepsilon\})=0.
\end{equation}
\end{df}

\begin{uw}
\label{uw1.1} Let $u,u_n,n\ge1$, be quasi-continuous. Let us
consider the following condition: for every $\varepsilon>0$,
\begin{equation}
\label{eq1.15} \lim_{n\rightarrow\infty}P_x(\sup_{t\ge0}
|u_n(X_t)-u(X_t)|>\varepsilon)=0
\end{equation}
for $m$-a.e. $x\in E$. Condition (\ref{eq1.15}) is equivalent to
the quasi-uniform, up to a subsequence, convergence of $\{u_n\}$
to $u$. To see this, let us set
$A_n^\varepsilon=\{|u_n-u|>\varepsilon\}$ and for arbitrary nearly
Borel set $B\subset E$ put $p_B(x)=P_x(\sigma_B<\infty)$, $x\in
E$. Assume that (\ref{eq1.15}) holds. By the diagonal method there
exists a subsequence (still denoted by $(n)$) such that
$p_{B^\varepsilon_n}(x)\rightarrow 0$, $m$-a.e. for every
$\varepsilon >0$, where $B^\varepsilon_n=\bigcup_{k\ge n}
A^\varepsilon_k$. Hence, by \cite[Corollary 4.3.1]{Fukushima},
CAP$(B^\varepsilon_n)\rightarrow 0$ for every $\varepsilon>0$,
which implies that $u_n\rightarrow u$ quasi-uniformly. Now assume
that $u_n\rightarrow u$ quasi-uniformly. Then by \cite[Theorem
2.1.5]{Fukushima}, $\EE(p_{A^\varepsilon_n},p_{A^\varepsilon_n})
=\mbox{CAP}(A^\varepsilon_n)\rightarrow0$. Therefore, up to a
subsequence, $p_{A^\varepsilon_n}\rightarrow 0$, $m$-a.e.  Let us
also mention that by the standard argument $``m\mbox{-a.e.}"$ in
condition (\ref{eq1.15}) may be replaced by ``q.e."
\end{uw}

\begin{uw}
Replacing CAP by Cap in (\ref{eq1.14}) we get a notion of
convergence which is weaker then the quasi-uniform convergence. In
fact, if
\begin{equation}
\label{eq4.3}
\lim_{n\rightarrow\infty}\mbox{Cap}(\{|u_n-u|>\varepsilon\})=0
\end{equation}
for every $\varepsilon>0$ then by \cite[Lemma IV.4.5]{MR},
$u_n\rightarrow u$ quasi-uniformly on every compact set $K\subset
E$. Therefore the convergence defined by (\ref{eq4.3}) may be
called a locally quasi-uniform convergence.
\end{uw}

\begin{stw}
\label{stw1.3} Let $\mu,\mu_n\in\mathcal{R}$ and let $u=R\mu$,
$u_n=R\mu_n$. If $u_n\rightarrow u$ quasi-uniformly then there
exists a subsequence (still denoted by $(n)$) such that for q.e.
$x\in E$,
\[
\lim_{n\rightarrow\infty}E_x\sup_{t\ge 0}
|A_t^{\mu_n}-A^{\mu}_t|=0.
\]
\end{stw}
\begin{dow}
Since $u_n(x)=E_x\int_0^\zeta\,dA_t^{\mu_n}\rightarrow u(x)$,
$\sup_{n\ge1}E_x\int_0^\zeta\,dA_t^{\mu_n}<\infty$, which when
combined with the quasi-uniform convergence of $\{u_n\}$ implies
that $\{u_n(X)\}$ satisfies the condition UT under $P_x$ for q.e.
$x\in E$ (see \cite[Proposition 3.2]{JMP}). Therefore by
\cite[Theorem 1.8]{Jacod} (see also \cite[Corollary 2.8]{JMP}),
for every $\varepsilon>0$,
\[
\lim_{n\rightarrow\infty}P_x(\sup_{t\ge 0}
|A_t^{\mu_n}-A_t^\mu|>\varepsilon)=0
\]
for q.e. $x\in E$. This and the fact that $u_n\rightarrow u$,
$m$-a.e. implies that the family $\{A^{\mu_n}_\zeta\}$ is
uniformly integrable with respect to $P_x$ for $m$-a.e. $x\in E$.
Applying  the Vitali theorem yields the
desired result.
\end{dow}

\begin{lm}
\label{lm1.1} Assume that $\mu,\mu_n\in S_0^{(0)}$ and
$\mu_n\rightarrow\mu$ strongly in $S_0^{(0)}$. Let $\{u_n\}$ be a
sequence of quasi-continuous functions such that $0\le u_n\le c$
for some $c>0$ and $u_n\rightarrow u$ quasi-uniformly. Then for
every positive $\eta\in L^2(E;m)$ and every $\alpha>0$,
\begin{equation}
\label{eq3.21} \int_E u_nR_\alpha\eta\,d\mu_n\rightarrow\int_E
uR_\alpha\eta\,d\mu.
\end{equation}
\end{lm}
\begin{dow}
Since $\mu_n\rightarrow\mu$ in $S_0^{(0)}$, it is easy to see that
$R\mu_n\rightarrow R\mu$ in the $\EE$-norm. Therefore by
\cite[Lemma 5.1.1]{Fukushima} there exists a subsequence (still
denoted by $(n)$) such that $R\mu_n\rightarrow R\mu$
quasi-uniformly. By this and Proposition \ref{stw1.3},
$E_x\sup_{t\ge 0}|A_t^{\mu_n}-A_t^\mu|\rightarrow0$ for q.e. $x\in
E$. Consequently,
\[
E_x\int_0^\zeta e^{-\alpha t}u_n(X_t)\,dA_t^{\mu_n}\rightarrow
E_x\int_0^\zeta e^{-\alpha t}u(X_t)\,dA_t^{\mu}
\]
for q.e. $x\in E$, so (\ref{eq3.21}) follows by the Lebesgue
dominated convergence theorem.
\end{dow}
\medskip

The following proposition is a generalization of
\cite[Theorem 1]{LiskevichSemenov}.
\begin{stw}
\label{stw1.2} Let  $\mu\in\mathcal{M}_{0,b}^+$ and $u=R\mu$. If
$\int_Eu^{p-1}d\mu<\infty$ for some $p>1$ then $u^{p/2}\in
D_e[\EE]$ and there exists $c_p>0$ such that
\[
\EE(u^{p/2},u^{p/2})\le c_p(u^{p-1},\mu).
\]
\end{stw}
\begin{dow}
Let $\theta\in D(A)$ be such that $0\le \theta\le 1$ and
$\theta\in L^1(E;m)$. Let us choose a nest $\{F_n\}$ such that
$\mathbf{1}_{F_n}\cdot\mu$, $\mathbf{1}_{F_n}u^{p-1}\cdot\mu\in
S^{(0)}_{00}$, $n\ge 1$, and by $u_n(\cdot;\lambda,\theta,\alpha)$
denote a solution of
\[
-A_\lambda u_n(\lambda,\theta,\alpha)=\theta\alpha R_\alpha\mu_n
\]
with $\mu_n=\mathbf{1}_{F_n}\cdot\mu$, $\alpha>0$ and
$A_\lambda=A-\lambda I$, $\lambda>0$.  Observe that $\theta\alpha
R_\alpha\mu_n\in L^2(E;m)\cap L^\infty(E;m)$. By \cite[Theorem
1]{LiskevichSemenov}, $u_n^{p/2}(\lambda,\theta,\alpha)\in D[\EE]$
and there exists $c_p>0$ such that
\begin{equation}
\label{eq1.11}
\EE(u_n^{p/2}(\lambda,\theta,\alpha),u_n^{p/2}(\lambda,\theta,\alpha))
\le c_p(u_n^{p-1}(\lambda,\theta,\alpha),\theta\alpha
R_\alpha\mu_n).
\end{equation}
Let $u_n(\cdot; \lambda,\theta)$ be a solution of
\[
-A_\lambda u_n(\lambda,\theta)=\theta\cdot\mu_n.
\]
By the very definition of  a solution,
\[
u_n(x;\lambda,\theta,\alpha)=E_x\int_0^\zeta e^{-\lambda r}
\Big(E_{X_r}\int_0^\zeta \alpha e^{-\alpha
t}\,dA_t^{\mu_n}\Big)\theta(X_r)\,dr
\]
and
\begin{equation}
\label{eq1.11a} u_n(x;\lambda,\theta)=E_x\int_0^\zeta e^{-\lambda
t}\theta(X_t)\,dA^{\mu_n}_t =E_x\int_0^\zeta e^{-\lambda
t}\theta(X_t)\mathbf{1}_{F_n}(X_t)\,dA^{\mu}_t
\end{equation}
for q.e. $x\in E$. Therefore by the Markov property and  Fubini's
theorem,
\begin{align*}
u_n(x;\lambda,\theta,\alpha)&=E_x\int_0^\zeta e^{-\lambda r}
\Big(E_x\int_r^\zeta\alpha
e^{-\alpha(t-r)}\,dA_t^{\mu_n}\Big)\theta(X_r)\,dr\\
&=E_x\int_0^\zeta \alpha e^{-\alpha t} \Big(\int_0^t
e^{(\alpha-\lambda)r}\theta(X_r)\,dr\Big)dA_t^{\mu}.
\end{align*}
Since $\theta\in D[\EE]$, $t\mapsto\theta(X_t)$ is c\`adl\`ag.
Therefore by standard calculations,
\[
\lim_{\alpha\rightarrow\infty}\alpha e^{-\alpha t}\int_0^t
e^{(\alpha-\lambda)r}\theta(X_r)\,dr
=\lim_{\alpha\rightarrow\infty}\int_0^t \alpha
e^{-\alpha(t-r)}e^{-\lambda r}\theta (X_r)\,dr= e^{-\lambda t}
\theta(X_t)
\]
and
\[
\alpha e^{-\alpha t}\int_0^t e^{(\alpha-\lambda)r}\theta(X_r)\,dr
\le 2 e^{-\lambda t}
\]
for $\alpha\ge\lambda$. Therefore applying the Lebesgue dominated
convergence theorem we get
\[
\lim_{\alpha\rightarrow\infty}u_n(x;\lambda,\theta,\alpha)
=\lim_{\alpha\rightarrow\infty} u_n(x;\lambda,\theta)
\]
for q.e $x\in E$. Observe that
\begin{equation}
\label{eq3.311} \|u_n^{p-1}(\lambda,\theta,\alpha)\|_\infty\le
\|R\mu_n\|^{p-1}_\infty:= c(n).
\end{equation}
Indeed, we have
\begin{align*}
u_n(x;\lambda,\theta,\alpha)\le R_\lambda(\alpha
R_\alpha(\mu_n))= \alpha
R_\alpha(R_\lambda(\mu_n))\le\alpha R_\alpha(\|R_\lambda
\mu_n\|_\infty)\le\|R\mu_n\|_\infty.
\end{align*}
From this and (\ref{eq1.11}) it follows that
\begin{align}
\label{eq1.12} \nonumber
\EE(u_n^{p/2}(\lambda,\theta,\alpha),u_n^{p/2}(\lambda,\theta,\alpha))
&\le c_p(u_n^{p-1}(\lambda,\theta,\alpha),\theta\alpha
R_\alpha(\mu_n))\\
& =c_p(\alpha R_\alpha
(u_n^{p-1}(\lambda,\theta,\alpha)\cdot\theta),\mu)\le
c_pc(n)\|\mu\|_{TV}
\end{align}
and
\begin{align}
\label{eq3.312}
\nonumber\EE(u_n(\lambda,\theta,\alpha),u_n(\lambda,\theta,\alpha))
&\le \EE_\lambda(u_n(\lambda,\theta,\alpha),u_n(\lambda,\theta,\alpha))\\
&\le (u_n(\lambda,\theta,\alpha),\alpha R_\alpha \mu_n)\le
c(n)^{1/(p-1)}\|\mu_n\|_{TV}.
\end{align}
Let us fix a sequence $\{\alpha_k\}\subset(0,\infty)$ such that
$\alpha_k\nearrow\infty$ and set
\[
S_k(u_n(x;\lambda,\theta,\alpha_k)) =\frac{1}{k}
\sum_{i=1}^ku_n(x;\lambda,\theta,\alpha_i).
\]
By (\ref{eq3.312}) and Mazur's theorem we may assume that
$S_k(u_n(\cdot;\lambda,\theta,\alpha_k))\rightarrow
u_n(\cdot;\lambda,\theta)$ in $\EE$. Therefore by \cite[Lemma
5.1.1]{Fukushima} and Remark \ref{uw1.1} there exists a
subsequence (still denoted by $(k)$) such that
$S_k(u_n(\lambda,\theta,\alpha_k))\rightarrow u_n(\lambda,\theta)$
quasi-uniformly as $k\rightarrow \infty$. It is an elementary
check that $\alpha_kR_{\alpha_k}(\mu_n)\rightarrow\mu_n$ weakly in
$S_0^{(0)}$ as $k\rightarrow \infty$. So, again by  Mazur's
theorem we may assume that $S_k(\alpha_k
R_{\alpha_k}\mu_n)\rightarrow \mu_n$ strongly in $S^{(0)}_0$.
Therefore by Lemma \ref{lm1.1}, up to a subsequence we have
\begin{equation}
\label{eq9.99}
(S_k^{p-1}(u_n(\lambda,\theta,\alpha_k))
\cdot\theta,S_k(\alpha_kR_{\alpha_k}\mu_n))
\rightarrow(u_n^{p-1}(\lambda,\theta)\cdot\theta,\mu_n)
\end{equation}
as $k\rightarrow\infty$. By
\cite[Theorem 1]{LiskevichSemenov},
\begin{align}
\label{eq1.1300} &\EE(S_k^{p/2}
(u_n(\lambda,\theta,\alpha_k)),S_k^{p/2}(u_n(\lambda,\theta,\alpha_k)))
\nonumber\\
&\qquad\le c_p(S_k^{p-1}(u_n(\lambda,\theta,\alpha_k)),
S_k(\alpha_kR_{\alpha_k}(\mu_n))\cdot\theta).
\end{align}
From this and (\ref{eq9.99}) we conclude that $\sup_{n\ge
1}\|S_k^{p/2}(u_n(\lambda,\theta,\alpha_k))\|_{\EE}<\infty$, which
implies that, up to subsequence,
$\{S_k^{p/2}(u_n(\lambda,\theta,\alpha_k))\}$ is weakly convergent
in $\EE$ to some $v\in D_e[\EE]$. Since by \cite[Lemma
5.1.1]{Fukushima} and Remark \ref{uw1.1} strong, up to a
subsequence, convergence  in $\EE$ implies quasi-uniform
convergence, by standard reasoning we get
$v=u^{p/2}_n(\lambda,\theta)$. Therefore by  (\ref{eq9.99}) and
\cite[Theorem 1]{LiskevichSemenov},
\begin{align}
\label{eq1.13}
\EE(u_n^{p/2}(\lambda,\theta),u_n^{p/2}(\lambda,\theta))
&\le\liminf_{k\rightarrow\infty}\EE(S_k^{p/2}
(u_n(\lambda,\theta,\alpha_k)),S_k^{p/2}(u_n(\lambda,\theta,\alpha_k)))
\nonumber\\
&\le
c_p\liminf_{k\rightarrow\infty}(S_k^{p-1}(u_n(\lambda,\theta,\alpha_k)),
S_k(\alpha_kR_{\alpha_k}(\mu_n))\cdot\theta)\nonumber \\
&=c_p(u_n^{p-1}(\lambda,\theta)\cdot\theta,\mu_n)
\le c_p(u_n^{p-1}(\lambda,\theta)\cdot\theta,\mu).
\end{align}
Let us choose $\theta_l\in D(A)$ such that $0\le\theta_l\le1$ and
$\theta_l\nearrow 1$.  For instance, one can take $\theta_l=l R_l
e_{F_l}$, where $e_{F_l}$ is the equilibrium function for the set
$F_l$ (see \cite[Chapter 2]{Fukushima}) and $\{F_l\}$ is defined
at the beginning of the proof. From  (\ref{eq1.11a})  and the fact that
\[
u(x)=E_x\int_0^\zeta dA^\mu_t,\quad x\in E
\]
one can deduce that
\[
u_n(x;\lambda,\theta_l)\le u,\qquad
\lim_{l\rightarrow\infty}\lim_{\lambda\rightarrow0}\lim_{n\rightarrow\infty}
u_n(x;\lambda,\theta_l)=u(x)
\]
for q.e. $x\in E$. This when combined with (\ref{eq1.13}) and the
assumptions of the proposition gives the desired result.
\end{dow}

\begin{tw}
\label{stw1.4} Assume that $u$ is a solution of
\mbox{\rm(\ref{eq1.1})}.
\begin{enumerate}
\item[\rm(i)] If $\mu\in S_{00}^{(0)}$ then $u\in L^\infty(E;m)$ and
\[
\|u\|_\infty\le c_2(\gamma+1)^{1/(\gamma+1)}
\|R\mu\|_\infty^{1/(\gamma+1)}.
\]
\item[\rm(ii)] If $\mu\in\MM_{0,b}^+(E)$ then
$u^{(\gamma+1)/2}\in D_e[\EE]$ and
\[
\|u^{(\gamma+1)/2}\|^2_\EE\le c(\gamma)c_2\|\mu\|_{TV.}
\]
\end{enumerate}
\end{tw}
\begin{dow}
(i) By the very definition of the space $S_{00}^{(0)}$, $R\mu\in
L^\infty(E;m)$. By the Meyer-Tanaka formula and (\ref{eq1.2}),
\[
u^{\gamma+1}(x)\le (\gamma+1)E_x\int_0^\zeta
g(u)u^\gamma(X_t)\,dA^{\mu}_t \le c_2(\gamma+1)R\mu(x),
\]
from which the desired estimate immediately follows. \\
(ii)  Let us put $\nu=g(u)\cdot\mu$ and $p=1+\gamma$. Then $p>1$
and
\[\int_Eu^{p-1}\,d\nu\le c_2\int_Eu^{p-1}\cdot\frac{1}{u^{p-1}}\,d\mu
=c_2\|\mu\|_{TV}.
\]
By the above estimate and Proposition \ref{stw1.2},
$u^{(\gamma+1)/2}\in D_e[\EE]$ and there exists $c(\gamma)>0$ such
that
\[
\EE(u^{(\gamma+1)/2},u^{(\gamma+1)/2})\le
c(\gamma)\int_Eu^\gamma\,d\nu\le c_2c(\gamma)\cdot\|\mu\|_{TV},
\]
which completes the proof.
\end{dow}

\begin{prz}
\label{ex4.7} Let $(\EE,D[\EE])$ be the form defined by
(\ref{eq4.dfp}) with $A=\Delta^{\alpha/2}$ for some
$\alpha\in(0,2]$ and bounded domain
$D\subset \BR^d$.\\
(i) We first give two examples of $\mu\in S^{(0)}_{00}$. Assume
that $d>\alpha$. By \cite[Example 2.2.1]{Fukushima}, if $\mu\in S$
and
\[
\sup_{x\in D}\int_D \frac{1}{|x-y|^{d-\alpha}}\,\mu(dy)<\infty
\]
then $\mu\in S^{(0)}_{00}$. For instance, if $\mu=f\cdot m$ and
$f\in L^p(D;m)$ with $p>d/{\alpha}$ then $f\cdot m\in
S^{(0)}_{00}$. Now, let $\alpha=2$ and let $\mu$ denote the
Riemannian volume measure on some $(d-1)$ dimensional submanifold
$\Sigma$ of $D$. Then extending $\mu$ by zero to the whole $D$ we
get $\mu\in S^{(0)}_{00}$ (see \cite[Example iv]{AJM}). In both
cases, by Theorem \ref{stw1.4}(i),  the solution $u$ of
(\ref{eq1.1}) bounded.
\smallskip\\
(ii) Let $\alpha\in(0,2\wedge d)$. By \cite[Theorem 4.4.4, Eq.
(1.5.19)]{Fukushima}, $D_e[\EE]=H^{\alpha/2}_0(D)$, where
\[
H^{\alpha/2}_0(D)=\{u\in L^2(\BR^d;dx); u=0 \mbox{  on  }
\BR^d\setminus D\ \mbox{  and  }
\int_{\BR^d}|\hat{u}(x)||x|^\alpha\,dx<\infty\}
\]
and $\hat{u}$ denotes the Fourier transform of $u$. Therefore if
$\mu\in\MM^+_{0,b}$ and $u$ is a solution of (\ref{eq1.1}) then
$u^{(\gamma+1)/2}\in H^{\alpha/2}_0(D)$.
\end{prz}

\nsubsection{Stability: General results I}
\label{sec5}

In Sections \ref{sec5}--\ref{sec7} we study stability of solutions
of the problem
\begin{equation}
\label{eq1.17} -Au_n=g(u_n)\cdot\mu_n,\quad u_n>0
\end{equation}
under different assumptions on the convergence of measures $\mu_n$
and the limit measure $\mu$. It is known that each
measure $\mu\in \MM_{b}$ admits a unique decomposition of the form
\begin{equation}
\label{eq5.2} \mu=\mu_d+\mu_c,
\end{equation}
where $\mu_d\in \MM_{0,b},\mu_c\in\MM_b$ and $\mu_c\bot$Cap. The
measure $\mu_d$ is called the diffuse part of $\mu$, whereas
$\mu_c$ the concentrated part of $\mu$.

In the present section we prove some general results on stability
in case $\mu_c\neq0$. Then in Section \ref{sec6} we investigate
the case where $\mu$ is smooth, i.e. $\mu_c=0$. Finally, in
Section \ref{sec7} we turn back to the case  $\mu_c\neq0$ but we
assume that $\mu_n$ are of the form $\mu_n=j_{1/n}\ast\mu$, where
$j$ is some mollifier, and that $A$ corresponds to some form $\EE$
on $L^2(D;dx)$ with $D\subset\BR^d$.

\begin{lm}
\label{lm4.sub} Let $\{u_n\}$ be a sequence of excessive functions
on $E$ such  that $u_n\rightarrow 0,\, m$-a.e. Then there exists a
subsequence $(n')\subset (n)$ such that $u_{n'}\rightarrow 0$ q.e.
\end{lm}
\begin{dow}
Without loss of generality we may assume that $u_n\le 1$, $n\ge1$.
Let $(n')\subset (n)$ be such that $\sum_{n'\ge 1}\int_E
u_{n'}\,d\pi<\infty$ (for the definition of $\pi$ see Section
\ref{sec2}). Let $E\setminus B=\{x\in E;\, u_{n'}(x)\rightarrow
0\}$ and let $F$ be a compact subset of $E$ such that $K\subset
B$. Then
\[
P_\pi(D_F<\zeta)\le P_\pi(\limsup_{n'\rightarrow\infty}
u_{n'}(X_{D_F})>0)=0.
\]
Indeed, since $u_n$ is an excessive function,
\[
P_\pi(u_{n'}(X_{D_F})>\varepsilon) \le\varepsilon^{-1}E_\pi
u_{n'}(X_{D_F})\le\varepsilon\int_E u_{n'}\,d\pi.
\]
Therefore $u_{n'}(X_{D_F})\rightarrow0$, $P_\pi$-a.e. by the
Borel-Cantelli lemma. Hence Cap$(F)=0$ by \cite[Theorem
IV.5.28]{MR}. Since $F\subset B$ was arbitrary, Cap$(B)=0$.
\end{dow}
\medskip

Let us recall that a sequence $\{\mu_n\}$ of Radon measures on $E$
converges to some Radon measure on $E$ in the narrow topology if
$\int_Ef\,\mu_n(dx)\rightarrow\int_E f\,\mu(dx)$ for every bounded
continuous $f:E\rightarrow\BR$. If the last convergence holds true
for every continuous  $f$ having compact support then we say that
$\{\mu\}$ converges to $\mu$ in the vague topology.

\begin{stw}
\label{stw1.5}  Assume that $(\EE,D[\EE])$ satisfies Meyer's
hypothesis {\rm{(L)}} and $g$ satisfies \mbox{\rm(\ref{eq1.2})}
with $\gamma=1$. Let  $\mu\in\MM^+_b$ be  such that
$\mu\bot\mbox{\rm Cap}$ and let $\{\mu_n\}\subset\MM^+_{0,b}$ be a
sequence such that $\sup_{n\ge1}\|\mu_n\|_{TV}<\infty$ and
$\mu_n\rightarrow\mu$ in the narrow topology. If $u_n$ is a
solution of the problem \mbox{\rm(\ref{eq1.17})} then there exists
a subsequence (still denoted by $(u_n)$) such that $u_n\rightarrow 0$
q.e.
\end{stw}
\begin{dow}
Let $\varepsilon>0$. Since $(\EE,D[\EE])$ is regular, there exists
$\psi_\varepsilon\in D[\EE]\cap C_c(E)$ such that
\begin{equation}
\label{eq1.16}
0\le\psi_\varepsilon\le 1,\quad 0
\le\int_E(1-\psi_\varepsilon)\,d\mu\le\varepsilon,\quad
\EE(\psi_\varepsilon,\psi_\varepsilon)\le\varepsilon
\end{equation}
(see \cite[Lemma 2.2.7]{Fukushima}). By Theorem \ref{stw1.4} and
\cite[Propositions 2.4, 2.11]{Kl:MA1}, $u_n\in D_e[\EE]$ and there
exists a subsequence (still denoted by $(n)$) such that $\{u_n\}$ is
convergent $m$-a.e. and weakly in $D_e[\EE]$ to some function
$u\in D_e[\EE]$. Since $u_n\in D_e[\EE]$,
\[
\EE(u_n,\eta)=(\eta,g(u_n)\cdot\mu_n)
\]
for every $\eta\in D[\EE]\cap\BB^+(E)$. For arbitrary but fixed
$k>0$ set
\[
\eta=T_k(u_n)(1-\psi_\varepsilon) =T_k(u_n)-\psi_\varepsilon
T_k(u_n)\in D_e[\EE],
\]
where $T_k$ is the truncature operator, i.e. $T_k(y)=((-k)\vee
y)\wedge k$, $y\in\BR$. Then
\begin{align*}
\EE(u_n,T_k(u_n)(1-\psi_\varepsilon))
&=(T_k(u_n)(1-\psi_\varepsilon),g(u_n)\cdot\mu_n)\\
&\le c_2\int_E\frac{T_k(u_n)(1-\psi_\varepsilon)}{u_n}\,d\mu_n\le
c_2\int_E(1-\psi_\varepsilon)\,d\mu_n.
\end{align*}
Also
\[
\EE(u_n,T_k(u_n)(1-\psi_\varepsilon))=\EE(u_n,T_k(u_n))-
\EE(u_n,T_k(u_n)\psi_\varepsilon).
\]
Since $(\EE,D[\EE])$ is a Dirichlet form, it is Markovian. Hence
\[
\EE(T_k(u_n),T_k(u_n))\le\EE(u_n,T_k(u_n))
\]
for $n\ge1$ and consequently,
\[
\EE(T_k(u_n),T_k(u_n))\le c_2\int_E(1-\psi_\varepsilon)\,d\mu_n
+\EE(u_n,T_k(u_n)\psi_\varepsilon).
\]
Since $u_n$ is a potential,
\[
\EE(u_n,T_k(u_n)\psi_\varepsilon)\le k\EE(u_n,\psi_\varepsilon).
\]
Therefore
\[
\EE(u_n,T_k(u_n)\psi_\varepsilon)\le k\sqrt{\EE(u_n,u_n)}\cdot
\sqrt{\EE(\psi_\varepsilon,\psi_\varepsilon)}\le
k\varepsilon^{1/2}\sqrt{\EE(u_n,u_n)}.
\]
By Theorem \ref{stw1.4}, $c:=\sup_{n\ge1}
\sqrt{\EE(u_n,u_n)}<\infty$. Hence
\[
\EE(T_k(u_n),T_k(u_n)) \le c_2\int_E(1-\psi_\varepsilon)\,d\mu_n+k
c\varepsilon^{1/2}.
\]
Letting $n\rightarrow \infty$ in the above inequality and using
(\ref{eq1.16}) we obtain
\[
\EE(T_k(u),T_k(u))\le c_2\varepsilon+kc\varepsilon^{1/2}.
\]
Since $k,\varepsilon >0$ were arbitrary, $u\equiv0$. The result
now follows from Lemma \ref{lm4.sub}.
\end{dow}
\medskip

Let $T\in \mathcal{T}$ and $\Lambda\in\FF_T$. Write
\[
T_\Lambda(\omega)=\left\{
\begin{array}{l}T(\omega),\quad \omega\in\Lambda, \\
\infty,\,\qquad \omega\notin \Lambda.
\end{array}
\right.
\]
It is well known (see \cite[Section III.2]{Protter}) that
$T_\Lambda\in\mathcal{T}$.

\begin{lm}
\label{lm.van} If $u\in D_e[\EE]$ then $u\in \mathbf{D}$ and for
every $\{T_n\}\subset \mathcal{T}$ such that $T_n\nearrow T\ge
\zeta$,
\[
u(X_{T_n})\rightarrow 0,\quad P_x\mbox{-a.s.}
\]
for q.e. $x\in E$.
\end{lm}
\begin{dow}
Let $\Lambda=\{\omega\in\Omega:
T_n(\omega)<\zeta(\omega),n\ge1\}$. Then $\Lambda\in\FF_T$ because
$\Lambda=\bigcap_{n\ge 1}\{T_n<\zeta\}$ and $\{T_n<\zeta\}\in
\FF_{T_n}\cap\FF_\zeta\subset\FF_T$. Also observe  that
$T=T_\Lambda\wedge T_{\Lambda^c}$ and that $T_\Lambda$ is
predictable. Since $u^+,u^-\in D[\EE]$, we may assume that
$u\ge0$. Let $v\in D_e[\EE]$ be an excessive function such that
$v\ge u$ q.e. (for the existence of such function see
\cite[Theorem I.2.6]{MR}). By \cite[Theorems 2.2.1,
5.1.1]{Fukushima} there exists a positive measure $\mu\in
S^{(0)}_{0}$ such that
\begin{equation}
\label{eq5.3} v(x)=E_x\int_0^\zeta dA^\mu_t
\end{equation}
for q.e. $x\in E$. By the strong Markov property,
\[
v(X_{T_n})=\int_{T_n\wedge\zeta}^\zeta dA^\mu_t
-\int_{T_n\wedge\zeta}^\zeta dM_t,\quad P_x\mbox{-a.s.}
\]
for q.e $x\in E$, where
\[
M_t=E_x\Big(\int_0^\zeta dA^\mu_r|\FF_t\Big)-v(X_0),\quad t\ge 0.
\]
By (\ref{eq5.3}), $v\in\mathbf{D}$. Consequently,  $u\in
\mathbf{D}$. Since $A^\mu$ is continuous,
\[
\int_{T_n\wedge\zeta}^\zeta dA^\mu_t\rightarrow 0, \quad
P_x\mbox{-a.s.}
\]
for q.e. $x\in E$. Moreover,
\[
\int_{T_n\wedge\zeta}^\zeta dM_r\rightarrow \Delta
M_\zeta\mathbf{1}_\Lambda =\Delta M_T\mathbf{1}_\Lambda=\Delta
M_{T_\Lambda}\mathbf{1}_\Lambda, \quad P_x\mbox{-a.s.}
\]
for q.e. $x\in E$. Since the filtration $(\FF_t)_{t\ge 0}$ is
quasi-left continuous, every  martingale with respect to it  has
only totally inaccessible jumps (see, e.g., \cite[Theorem
A.3.6]{Fukushima}). Hence $\Delta
M_{T_\Lambda}\mathbf{1}_\Lambda=0$, $P_x$-a.s. since $T_\Lambda$
is predictable. This proves the lemma.
\end{dow}
\medskip

The next general stability result will play an important role in
the proof of Theorem \ref{tw1.6}, which in turn is used in the
proof of our main Theorem \ref{tw1.7} on existence of solutions of
(\ref{eq1.1}) with general bounded Borel measure on the right-hand
side. Perhaps it is also appropriate to make here the following
general  comments.


In most  papers devoted to stability of solutions of semilinear
equations with measure data the following equation
\begin{equation}
\label{eq.additive} -\Delta u=f(x,u)+\mu
\end{equation}
is considered. Let $\{\mu_n\}$ be an approximation of a
nonnegative measure $\mu$ in the narrow topology and let $u_n$ be
a solution of (\ref{eq.additive}) with $\mu$ replaced by $\mu_n$.
Usually the limit $u$ of $\{u_n\}$ depends on the form of the
approximation of $\mu$ (see \cite{MP}). To be more precise, the
limit $u$ solves (\ref{eq.additive}) with $\mu$ replaced by some
nonnegative Borel measure $\mu^{\#}$, depending on $\{\mu_n\}$,
such that $\mu^{\#}\le\mu$ ($\mu^{\#}$ is called  the reduced
limit of $\{\mu_n\}$). The question naturally arises whether
similar phenomenon takes place in case of equations of the form
(\ref{eq1.1}). In \cite{BoccardoOrsina} it is observed that in the
particular case of equation (\ref{eq1.1}) with $A=\Delta$, $g$
satisfying (\ref{eq1.2}) with $\gamma\ge1$ and singular $\mu$
(i.e. $\mu=\mu_c)$ we have that $u_n\rightarrow 0$ for any
approximation of $\mu$ by uniformly bounded measures $\mu_n$ such
that $\mu_n\rightarrow\mu$ in the narrow topology. In different
words, for any approximation of $\mu$ in the limit equation the
whole singular part of $\mu$ disappear. We do not know whether
similar result holds true for any $\gamma>0$ and/or general
Dirichlet operator $A$. However, in Theorem \ref{tw1.5} below we
are able to prove a related result for general $A$ and bounded
measure. It says that the limit function $u$ satisfies an equation
with a measure $\nu$ on the right-hand side which is always smooth
independently on the approximation of $\mu$. But let us stress
that Theorem \ref{tw1.5} does not imply the result of
\cite{BoccardoOrsina}, because even in case $\mu=\mu_c$ we do not
know whether $\nu=0$. It is also worth mentioning that in Theorem
\ref{tw1.5} we consider the convergence in the vague topology.

In the proof of Theorem \ref{tw1.5} we will need the following
additional notation. For every open set $U\subset E$ we write
\[
D_{e,U}[\EE]=\{u\in D_e[\EE]:u=0\mbox{ q.e. on } E\setminus U\}.
\]
It is known (see \cite[Theorems 4.4.3, 4.4.4]{Fukushima}) that the
pair $(\EE,D_{e,U}[\EE])$ is again a regular transient symmetric
Dirichlet form. By $\{R^U_\alpha,\alpha\ge 0\}$ we denote the
resolvent associated with $(\EE,D_{e,U}[\EE])$. For a compact set
$F\subset U$ we denote by $e^U_F$ the equilibrium function
associated with $(\EE,D_{e,U}[\EE])$ and  $F$. By \cite[Theorem
2.1.5]{Fukushima}, $e^U_F$ is quasi-continuous and
\[
e^U_F=1\mbox{  q.e.  on  } F,\quad 0\le e^U_F\le 1\mbox{
q.e.},\quad e^U_F\in D_{e,U}[\EE]\subset D_e[\EE].
\]
The last property implies that $e^U_F=0$ q.e. on $E\setminus U$.
\begin{tw}
\label{tw1.5} Assume that $(\EE,D[\EE])$ satisfies Meyer's
hypothesis $\mbox{\rm{(L)}}$. Let $\mu\in\MM_b^+$,
$\{\mu_n\}\subset\MM_{0,b}^+$ be such that
$\sup_{n\ge1}\|\mu_n\|_{TV}<\infty$ and $\mu_n\rightarrow\mu$
vaguely. Let $u_n$ be a solution of \mbox{\rm(\ref{eq1.17})} and
let $\nu_n=g(u_n)\cdot\mu_n$. Then
\begin{enumerate}
\item[\rm(i)]$\{\nu_n\}$ is tight in the vague topology and its each  limit
point $\nu$ belongs to $\mathcal{R}$,
\item[\rm(ii)] if $\nu_{n'}\rightarrow \nu$ vaguely for some subsequence
$(n')\subset(n)$ then there is a further subsequence $(n'')\subset
(n')$ such that $u_{n''}\rightarrow u$, $m$-a.e., where $u$ is a
solution of
\[
-Au=\nu.
\]
\end{enumerate}
\end{tw}
\begin{dow}
Since $u_n$ is a solution of (\ref{eq1.17}), it is
quasi-continuous, $u_n\in \mathbf{D}$  and by the Markov property
there is a martingale additive functional $M^n$ of $\BX$ such that
\[
u_n(X_t)=\int_t^\zeta g(u_n)(X_r)\,dA_r^{\mu_n}-
\int_t^\zeta\,dM_r^n,\quad t\in[0,\zeta],\quad P_x\mbox{-a.s.}
\]
for q.e. $x\in E$. By the Meyer-Tanaka formula,
\begin{align*}
u_n^{\gamma+1}(X_t)+\int_t^\zeta\,dK_r^\gamma&=(\gamma+1)
\int_t^\zeta u_n^\gamma\cdot g(u_n)(X_r)\,dA_r^{\mu_n}\\
&\quad -(\gamma+1)\int_t^\zeta u_n^\gamma(X_{r-})\,dM_r^n,\quad
t\in[0,\zeta],\quad P_x\mbox{-a.s.}
\end{align*}
for some increasing process $K^\gamma$ such that $K_0^\gamma=0$.
Therefore by (\ref{eq1.2}),
\begin{equation}
\label{eq1.18} u_n^{\gamma+1}(X_t)\le c_2(\gamma+1)
E_x\Big(\int_0^\zeta\,dA_r^{\mu_n}|\FF_t\Big), \quad
t\in[0,\zeta],\quad P_x\mbox{-a.s.}
\end{equation}
for q.e. $x\in E$. In particular, for every $\beta\in
S_{00}^{(0)}$,
\begin{equation}
\label{eq1.19} \int_E u_n^{(\gamma+1)}(x)\,d\beta(x)\le
c_2(\gamma+1)\|R\beta\|_\infty\|\mu_n\|_{TV}.
\end{equation}
Observe that
\[
u_n(x)=E_x\int_0^\zeta\,dA_t^{\nu_n}
\]
for q.e. $x\in E$. Therefore from (\ref{eq1.19}) it follows
that for every   $\beta\in S_{00}^{(0)}$,
\begin{equation}
\label{eq1.25}
\sup_{n\ge1}\int_ER\beta\,d\nu_n=\sup_{n\ge1}\int_Eu_n\,d\beta<\infty.
\end{equation}
Let $K\subset E$ be a compact set. By \cite[Lemma
2.2.6]{Fukushima}, $e_K=R\beta_K$ for some $\beta_K\in
S^{(0)}_{00}$, where $e_K$ is the equilibrium function for $K$.
Since $e_K$ is positive and $e_K(x)=1$ q.e. on $K$, we conclude
from (\ref{eq1.25})  that $\{\nu_n\}$ is tight in the vague
topology. Let $\nu$ denote a limit point of $\{\nu_n\}$. By
(\ref{eq1.18}) and \cite[Lemma 6.1]{BDHPS}, for every $q\in (0,1)$
we have
\[
E_x\sup_{t\ge 0} u_n^{q(\gamma+1)}(X_t)\le
c_2^q\frac{(\gamma+1)^q}{(1-q)}
\Big(E_x\int_0^\zeta\,dA_t^{\mu_n}\Big)^q
\]
for q.e. $x\in E$. It follows that for every $\beta\in
S_{00}^{(0)}$ such that $\beta(E)=1$,
\begin{align}
\label{eq1.20} E_\beta\sup_{t\ge 0} |u_n(X_t)|^\alpha\le
c_2^q\frac{(\gamma+1)^q}{1-q}
\Big(E_\beta\int_0^\zeta\,dA_t^{\mu_n}\Big)^q\le
c_2^q\frac{(\gamma+1)^q}{1-q}\|R\beta\|^q_\infty\|\mu_n\|^q_{TV},
\end{align}
where $\alpha=q(\gamma+1)$. By Theorem \ref{stw1.4},
\begin{equation}
\label{eq1.21} \sup_{n\ge 1}
\EE(u_n^{(\gamma+1)/2},u_n^{(\gamma+1)/2})<\infty.
\end{equation}
By \cite[Lemma 94, page 306]{DellacherieMeyer} there exists a
subsequence (still denoted by $(n)$) such that
\begin{equation}
\label{eq1.22} u_n\rightarrow u,\quad m\mbox{-a.e.},
\end{equation}
where $u$ is an excessive function. By (\ref{eq1.21}),
$u^{(\gamma+1)/2}\in D_e[\EE]$. Therefore by Lemma \ref{lm.van},
$u^{(\gamma+1)/2}(X_{T_n})\rightarrow 0$ for every sequence
$\{T_n\}$ of stopping times such that $T_n\nearrow T\ge\zeta$.
This implies that for q.e. $x\in E$,
\begin{equation}
\label{eq1.23}
u(X_{T_n})\rightarrow 0,\quad P_x\mbox{-a.s.}
\end{equation}
A key step in showing that $\nu$ is smooth is the proof that $u$
is a potential. We first prove the last property in the simpler
case where (\ref{eq1.2}) is satisfied for some $\gamma\ge 1$.
Since $u^{(\gamma+1)/2}\in D[\EE]$, it belongs to $\mathbf{D}$ by
Lemma \ref{lm.van}. Therefore by (\ref{eq1.23}),
\[
E_xu^{(\gamma+1)/2}(X_{T_n})\rightarrow0
\]
for q.e. $x\in E$. From this we conclude that if $\gamma\ge1$ then
for q.e. $x\in E$,
\begin{equation}
\label{eq1.24} (E_xu(X_{T_n}))^{(\gamma+1)/2}\le
E_xu^{(\gamma+1)/2}(X_{T_n})\rightarrow 0,
\end{equation}
so if (\ref{eq1.2}) with $\gamma \ge 1$ is satisfied then $u$ is a
potential. Now we turn to the case $\gamma\in(0,1)$. It is perhaps
worth explaining why it differs from the case $\gamma\ge1$. To
show that $u$ is a potential we have to know that
$E_xu(X_{T_n})\rightarrow 0$. This may be concluded from
(\ref{eq1.23}) if $u\in\mathbf{D}$. Unfortunately, the last
assertion cannot be concluded from the fact $u^{(\gamma+1)/2}\in
D[\EE]$ when (\ref{eq1.2}) is satisfied with $\gamma\in(0,1)$. Now
we give an alternative way to prove that $u\in\mathbf{D}$. It is
independent of the value of $\gamma>0$. For $x\in E$ write
\[
\lambda^\alpha_x=\delta_{\{x\}}\circ R_\alpha.
\]
Since $(\EE,D[\EE])$ satisfies Meyer's hypothesis (L),
$\lambda^\alpha_x\ll m$ for every $x\in E$ . Moreover, since $u_n$
is a quasi-continuous excessive function, $\alpha R_\alpha
u_n(x)\le u_n(x)$ for q.e $x\in E$. From this and (\ref{eq1.22})
it follows that for q.e. $x\in E$,
\begin{align}
\label{eq4.11} \liminf_{n\rightarrow\infty}
u_n(x)\ge\liminf_{n\rightarrow\infty}\alpha R_\alpha
u_n(x)&=\liminf_{n\rightarrow\infty}\int_E\alpha u_n(y)
\lambda^\alpha_x(y)\,m(dy) \nonumber\\ &\ge\int_E\alpha
u(y)\lambda^\alpha_x(y)\,m(dy)= \alpha R_\alpha u(x).
\end{align}
Since $u^{(\gamma+1)/2}\in D[\EE]$, $u$ is quasi-continuous. Hence
$\alpha R_\alpha u(x)\nearrow u(x)$ for q.e. $x\in E$ as
$\alpha\nearrow\infty$. Therefore (\ref{eq4.11}) implies that
\begin{equation}
\label{eq.dml} \liminf_{n\rightarrow\infty} u_n(x)\ge u(x)
\end{equation}
for q.e. $x\in E$. By the above, for q.e. $x\in E$ we have
\[
u^\alpha(X_t)\le\liminf_{n\rightarrow\infty} u_n^\alpha(X_t),\quad
t\ge 0,\quad P_x\mbox{-a.s.}
\]
Hence
\[
\sup_{t\ge 0}u^\alpha(X_t)\le\sup_{t\ge
0}\liminf_{n\rightarrow\infty} u_n^\alpha(X_t)\le
\liminf_{n\rightarrow\infty}\sup_{t\ge 0}u_n^\alpha(X_t).
\]
By  Fatou's lemma,
\[
E_\beta\sup_{t\ge 0}u^\alpha(X_t) \le\liminf_{n\rightarrow\infty}
E_\beta\sup_{t\ge 0}u_n^\alpha(X_t) \le\sup_{n\ge1}
E_\beta\sup_{t\ge 0}u_n^\alpha(X_t),
\]
so by (\ref{eq1.20}),
\[
E_\beta\sup_{t\ge 0}u^\alpha(X_t) \le
c_2^q\frac{(\gamma+1)^q}{1-q}\|R\beta\|^q_{\infty} (\sup_{n\ge
1}\|\mu_n\|^q_{TV})
\]
for every $\beta\in S_{00}^{(0)}$ such that $\beta(E)=1$. Since
$\alpha>1$, we get in particular that $u\in\mathbf{D}$.  Therefore
by (\ref{eq1.23}), for q.e. $x\in E$,
\[
E_x u(X_{T_n})\rightarrow 0
\]
for every $\{T_n\}\subset\mathcal{T}$ such that $T_n\nearrow
T\ge\zeta$, which implies that $u$ is a potential. Therefore by
\cite[Theorem IV.4.22]{BG} and \cite[Theorem 5.1.4]{Fukushima}
there exists a smooth measure $\bar{\nu}$ such that
\begin{equation}
\label{eq1.bar} u(x)=E_x\int_0^\zeta\,dA_t^{\bar{\nu}}
\end{equation}
for q.e. $x\in E$.  Let $\beta\in S^{(0)}_{00}$.  Then by
(\ref{eq1.19}), (\ref{eq.dml}) and \cite[Theorem 26, page
28]{DellacherieMeyer1} there exists $v\in L^1(E;\beta)$ such that
\[
\int_E|u_n-v|\,d\beta\rightarrow 0.
\]
By (\ref{eq1.21}), (\ref{eq1.22}) and the Banach-Saks theorem,
\[
\frac{u^\gamma_1+\ldots+u^\gamma_n}{n} \rightarrow
u^\gamma\quad\mbox{in   }(\EE,D_e[\EE])
\]
as $n\rightarrow\infty$. By this and \cite[Theorem
2.1.4]{Fukushima} we may assume that the above convergence holds
q.e. Hence
\[
\frac{u^\gamma_1+\ldots+u^\gamma_n}{n} \rightarrow u^\gamma,\quad
\beta\mbox{-a.e.}
\]
From this we easily deduce that $v=u$, $\beta$-a.e. Consequently,
\begin{equation}
\label{eq.rcon} (R\beta,\nu_n)\rightarrow( R\beta,\bar{\nu}).
\end{equation}
By Dynkin's formula (see \cite[(4.4.3)]{Fukushima}) and
\cite[Section 2.3]{Fukushima}, for every open $U\subset E$,
\[
R^U\beta=R\beta-R(\beta)_{E\setminus U},
\]
where $(\beta)_{E\setminus U}$ is the sweeping out of $\beta$ on
$E\setminus U$. Since $R(\beta)_{E\setminus U}\le R\beta$, we have that
$(\beta)_{E\setminus U}\in S^{(0)}_{00}$ because by \cite[Lemma
5.4]{KR:JFA}, $(\beta)_{E\setminus U}(E)\le \beta(E)$. Therefore from
(\ref{eq.rcon}) it follows that
\begin{equation}
\label{eq.rconu}
( R^U\beta,\nu_n)\rightarrow( R^U\beta,\bar{\nu}).
\end{equation}
Let $\Pi=\{F\subset E: F\mbox{-compact},\, \nu(\partial F)=0\}$.
Then $\Pi$ is a $\pi$-system and $\sigma(\Pi)=\BB(E)$. For
$F\in\Pi$ let $F_\varepsilon=\{x\in E;\,
\mbox{dist}(x,F)<\varepsilon\}$. Since $E$ is locally compact,
there exists $\varepsilon>0$ such that $F_\varepsilon$ is
relatively compact. By \cite[Lemma 2.2.6]{Fukushima} and comments
following it, $e^{F^\varepsilon}_F=R^{F_\varepsilon}\beta$ for
some  $\beta\in S^{(0)}_{00}$, so by (\ref{eq.rconu}) we
have
\[
\bar{\nu}(F_\varepsilon)\ge ( e^{F_\varepsilon}_F,\bar{\nu})\ge
\liminf_{n\rightarrow \infty}\nu_n(\mbox{Int}F)\ge
\nu(\mbox{Int}F)=\nu(F).
\]
Since $\varepsilon>0$ can be made  arbitrarily small, it follows
that $\bar{\nu}(F)\ge \nu(F)$ for $F\in\Pi$. On the other hand,
again by (\ref{eq.rconu}),
\[
\bar{\nu}(F)\le  ( e^{F_\varepsilon}_F,\bar{\nu}) \le
\limsup_{n\rightarrow 0}\nu_n(\overline{F_\varepsilon}) \le
\nu(\overline{F_\varepsilon}).
 \]
Hence $\bar{\nu}(F)\le \nu(F)$, $F\in\Pi$. Therefore
$\bar{\nu}(F)=\nu(F)$ for $F\in\Pi$, which implies that
$\bar{\nu}=\nu$.
\end{dow}

\nsubsection{Stability: General results II} \label{sec6}

In the further study of stability an important role will be played
by a new type of convergence of measures of the class $\RR$, which
we define below. Since this convergence is related to the uniform
convergence of associated additive functionals, we will denote it
by $\uAF$.

\begin{df}
Let $\mu_n,\mu\in\RR$. We say that $\mu_n\uAF\mu$ if for  every
sequence $(n')\subset (n)$ there exists a further subsequence
$(n'')\subset (n')$ such that
\[
\lim_{n''\rightarrow\infty}E_x\sup_{t\ge 0}
|A^{\mu_{n''}}_t-A^\mu_t|=0
\]
for q.e. $x\in E$.
\end{df}

\begin{stw}
\label{stw1.11} Assume that $(\EE,D[\EE])$ satisfies hypothesis
{\rm(L)}.  Let $\mu_n,\mu\in\RR$ and let $u_n,u$ be solutions of
\begin{equation}
\label{eq4.iff} -Au_n=\mu_n,\quad -Au=\mu.
\end{equation}
If $\mu_n\uAF\mu$ then there exists a subsequence (still denoted
by $(u_n)$) such that $u_n\rightarrow u$ quasi-uniformly.
\end{stw}
\begin{dow}
By the assumption, up to a subsequence we have
\begin{equation}
\label{eq1.46} E_x\sup_{t\ge 0}|A_t^{\mu_n}-A_t^\mu|\rightarrow0
\end{equation}
for q.e. $x\in E$. By (\ref{eq4.iff}) and the definition of a
solution,
\[
\sup_{t\ge 0}|u_n(X_t)-u(X_t)| \le \sup_{t\ge 0}E_x(\sup_{r\ge
0}|A_r^{\mu_n}-A_r^\mu||\FF_t).
\]
From this and \cite[Lemma 6.1]{BDHPS}, for every $q\in(0,1)$ we
have
\[
E_x\sup_{t\ge 0}|u_n(X_t)-u(X_t)|^q \le\frac{1}{1-q}(E_x\sup_{t\ge
0}|A_t^{\mu_n}-A_t^\mu|)^q.
\]
By (\ref{eq1.46}), for q.e. $x\in E$ the right-hand side of the
above inequality converges to zero as $n\rightarrow\infty$, which
by Remark \ref{uw1.1}  completes the proof.
\end{dow}

\begin{lm}
\label{lm3.cap11} Let $u$ be a quasi-continuous function. Then
\[
\lim_{k\rightarrow\infty}\mbox{\rm Cap}(\{u>k\})=0.
\]
\end{lm}
\begin{dow}
Let $\{F_n\}$ be a nest such that $u$ is continuous on $F_n$ for
every $n\ge1$. Since $(\EE,D[\EE])$ is a regular Dirichlet form,
the capacity Cap generated by $(\EE,D[\EE])$ is tight (see
\cite[Remark IV.3.2]{MR}), i.e. there exists a nest
$\{\tilde{F}_m\}$ of compact subsets of $E$ such that
Cap$(E\setminus \tilde{F}_m)\rightarrow 0$ as
$m\rightarrow\infty$. By subadditivity of the capacity Cap,
\[
\mbox{Cap}(u>k)\le \mbox{Cap}(E\setminus
F_n)+\mbox{Cap}(E\setminus \tilde{F}_m)
+\mbox{Cap}(F_n\cap\tilde{F}_m,u>k).
\]
Since $F_n\cap\tilde{F}_m$ is compact and $u$ is continuous on
$F_n$,  $u$ is bounded on $F_n\cap\tilde{F}_m$. Hence
$\mbox{Cap}(F_n\cap\tilde{F}_m,\,u>k)\rightarrow 0$ as
$k\rightarrow \infty$. The other two terms converge to zero by the
definition of the nest.
\end{dow}
\medskip

Theorem \ref{stw1.6} below will play a key role in the proof of
our main result on  existence of solutions of  (\ref{eq1.1}) with
general $\mu\in\MM_b$ (Theorem \ref{tw1.7}). It is worth pointing
out that Theorem \ref{stw1.6} is new even in case $A=\Delta$.

\begin{tw}
\label{stw1.6} Assume that $(\EE,D[\EE])$ satisfies $(\EE.7)$ and
Meyer's hypotheses {\rm{(L)}}. Let $\mu\in\MM^+_{0,b}$ be nontrivial and let
$\{\mu_n\}\subset\MM_{0,b}^+$ be a sequence such that
$\sup_{n\ge1}\|\mu_n\|_{TV}<\infty$ and $\mu_n\uAF\mu$. If $u_n,u$
are solutions  of
\begin{equation}
\label{eq1.27} -Au_n=g(u_n)\cdot\mu_n,\quad u_n>0,\quad
-Au=g(u)\cdot\mu,\quad u>0
\end{equation}
with $g$ satisfying {\rm{(H)}} then there exists a subsequence
(still denoted by $(u_n)$) such that $u_n\rightarrow u$ q.e.
\end{tw}

\begin{dow}
Let $g_1(u)=g(u)\wedge 1$, $u>0$, and let $w_n$ be a solution of
the problem
\[
-Aw_n=g_1(w_n)\mu_n,\quad w_n>0.
\]
By Proposition \ref{stw1.1}, $w_n\le u_n$ q.e. and $w_n\le v_n$
q.e., where
\[
v_n(x)=E_x\int_0^\zeta\,dA^{\mu_n}_t,\quad x\in E.
\]
Put
\[
v(x)=E_x\int_0^\zeta\,dA^{\mu}_t,\quad x\in E.
\]
By Proposition \ref{stw1.11}, up to a subsequence, $v_n\rightarrow
v$ quasi-uniformly. By the Meyer-Tanaka formula and (\ref{eq1.2}),
for $k\ge c_1^{1/\gamma}$ we have
\begin{align*}
w_n(x)\wedge k&\ge E_x\int_0^\zeta\mathbf{1}_{\{w_n\le k\}}(X_t)
\Big(\frac{c_1}{w_n^\gamma(X_t)}\wedge 1\Big)\,dA_t^{\mu_n}\\
&\ge E_x\int_0^\zeta\mathbf{1}_{\{w_n\le k\}}(X_t)
\Big(\frac{c_1}{k^\gamma}\wedge 1\Big)\,dA_t^{\mu_n} \ge
\frac{c_1}{k^\gamma}E_x\int_0^\zeta\mathbf{1}_{\{v_n\le
k\}}(X_t)\,dA_t^{\mu_n}.
\end{align*}
Let $\{F_m\}$ be a nest such that $v_n\rightarrow v$ uniformly on
$F_m$ for every $m\ge 1$. For $\varepsilon>0$ let us choose
$n(\varepsilon,m)$ so that $|v_n(x)-v(x)|\le\varepsilon$ for $x\in
F_m$ and $n\ge n(\varepsilon,m)$. Then
\[w_n(x)\wedge k\ge  \frac{c_1}{k^\gamma}E_x
\int_0^\zeta\mathbf{1}_{F_m}(X_t) \mathbf{1}_{\{v\le
k-\varepsilon\}}(X_t)\,dA_t^{\mu_n}.
\]
Let $\eta>0$ be such that $R\eta\le 1$. Write
$C^m_{k,\varepsilon}=F_m\cap\{v\le k-\varepsilon\}$ and
$\eta_{k,\varepsilon}^m(x)
=E_x\int_0^{\tau_{C^m_{k,\varepsilon}}}\eta(X_t)\,dt$. Then
$\eta_{k,\varepsilon}^m\in D_e[\EE]$ and
$\eta_{k,\varepsilon}^m\le 1$ q.e., $\eta_{k,\varepsilon}^m=0$
q.e. on $E\setminus C^m_{\varepsilon,k}$. Hence for $n\ge
n(\varepsilon,m)$,
\[
w_n(x)\wedge k\ge \frac{c_1}{k^\gamma}
E_x\int_0^\zeta\eta_{k,\varepsilon}^m(X_t)\,dA_t^{\mu_n}
\]
for q.e $x\in E$.
By \cite[Theorem IV.5.28]{MR}
$P_x(\lim_{k,m\rightarrow \infty}\tau_{C^m_{k,\varepsilon}}<\zeta)=0$
for q.e. $x\in E$.
Hence
\begin{equation}
\label{eq5.strp}
\eta^m_{k,\varepsilon}(x)\nearrow E_x\int_0^\zeta\eta(X_t)\,dt>0,
\quad\mbox{for q.e.  }x\in E.
\end{equation}
Let
\[
\psi_{n,\varepsilon}^m(x)
=E_x\int_0^\zeta\eta_{k,\varepsilon}^m(X_t)\,dA_t^{\mu_n},\quad
\psi_\varepsilon^m(x)
=E_x\int_0^\zeta\eta_{k,\varepsilon}^m(X_t)\,dA_t^{\mu}.
\]
Then for $n\ge n(\varepsilon,m)$,
\[
0\le \psi_{n,\varepsilon}^m(x) \frac{c_1}{k^\gamma} \le
u_n(x)\wedge k
\]
for q.e. $x\in E$. By the assumptions, up to a subsequence,
\begin{equation}
\label{eq4.17} \lim_{n\rightarrow\infty}E_x\sup_{t\ge
0}|\int_0^t\eta_{k,\varepsilon}^m(X_r)\,d(A_r^{\mu_n}-A_r^{\mu})
|=0
\end{equation}
for q.e. $x\in E$. By \cite[Lemma 6.1]{BDHPS},
\[
E_x\sup_{t\ge 0}|\psi_{n,\varepsilon}^m(X_t)
-\psi_{\varepsilon}^m(X_t)|^q\le\frac{1}{1-q}\Big(E_x\sup_{t\ge
0}|
\int_0^t\eta_{k,\varepsilon}^m(X_r)\,d(A_r^{\mu_n}-A_r^{\mu})|\Big)^q.
\]
This together with (\ref{eq4.17}) and Remark \ref{uw1.1} shows
that up to a subsequence,
\begin{equation}
\label{eq4.ab} \psi_{n,\varepsilon}^m\rightarrow\psi_\varepsilon^m
\quad\mbox{quasi-uniformly as  } n\rightarrow \infty.
\end{equation}
Therefore for every $\varepsilon >0,\, m\ge 1$ there exists a nest
$\{F^{\varepsilon,m}_j,\, j\ge 1\}$ such that
$\psi^m_{n,\varepsilon}\rightarrow \psi^m_\varepsilon$ uniformly
on $F^{\varepsilon,m}_j$ for every $j\ge 1$. By \cite[Lemma 94,
page 306]{DellacherieMeyer} there exists a subsequence (still
denoted by $(n)$) such that $\{u_n\}$ converges $m$-a.e. Now we
will show that  one can choose a subsequence such that $\{u_n\}$
converges q.e. To this end, for  $a>0$ set
\[
B^{n,m}_a=\{u_n\ge\frac 1 a\}\cap F_m,\qquad A_{a,\varepsilon}^m
=\{\frac{c_1}{k^\gamma}\psi^m_\varepsilon \ge\frac1a
-\varepsilon\}\cap F_m\cap F^{\varepsilon,m}_{j(\varepsilon,m)}
\]
and
\[
D^{m}_{a,\varepsilon}=A^m_{a,\varepsilon}\setminus(E\setminus
A^m_{a,\varepsilon})^r,
\]
 where $(E\setminus A^m_{a,\varepsilon})^r$ is the set of
regular points for $E\setminus A^m_{a,\varepsilon}$ (see
\cite{Fukushima}) and $j(\varepsilon,m)$ is such that
Cap$(E\setminus
F^{\varepsilon,m}_{j(\varepsilon,m)})<\varepsilon/m$.
It is known that $D^m_{a,\varepsilon}$
is the fine interior of $A^m_{a,\varepsilon}$. Then
\[
B^{n,m}_a\supset A^m_{a,\varepsilon},\quad n\ge n(m,\varepsilon)
\]
and $\{p_t^{D^{m}_{a,\varepsilon}},\,t\ge 0\}$ is a strongly
continuous semigroup on $L^2(D^{m}_{a,\varepsilon};m)$, where
\[
p_t^{D^{m}_{a,\varepsilon}}u(x)
=E_xu(X_t)\mathbf{1}_{\{t<\tau_{D^{m}_{a,\varepsilon}}\}}, \quad
x\in D^{m}_{a,\varepsilon}.
\]
By the probabilistic definition of a solution of (\ref{eq1.27}) we
have
\begin{align*}
|u_n(x)-E_xu_n(X_{t\wedge{\tau_{D^m_{a,\varepsilon}}}})|&=
E_x\int_0^{t\wedge\tau_{D^m_{a,\varepsilon}}}g(u_n)(X_r)\,dA_r^{\mu_n}\\
&\le c_2\int_0^{t\wedge\tau_{D^m_{a,\varepsilon}}}
\frac{1}{|u_n|^\gamma}(X_r)\,dA_r^{\mu_n} \le
c_2a^\gamma\int_0^{t\wedge\tau_{D^m_{a,\varepsilon}}}\,dA_r^{\mu_n}.
\end{align*}
Hence
\[
\lim_{t\rightarrow 0^+}\sup_{n\ge 1}
|u_n(x)-E_xu_n(X_{t\wedge\tau_{D^m_{a,\varepsilon}}})|\le c_2\cdot
a^\gamma\lim_{t\rightarrow 0^+}\sup_{n\ge
1}E_x\int_0^t\,dA_r^{\mu_n}.
\]
Since $\mu_n\uAF\mu$, for every $\delta>0$ there exists
$n_0\in\BN$ such that for every $t\ge 0$,
\[
|E_x\int_0^t\,dA_r^{\mu_n}-E_x\int_0^t\,dA_r^\mu|\le\delta,\quad
n\ge n_0.
\]
Therefore
\begin{align*}
\lim_{t\rightarrow 0^+}\sup_{n\ge 1}
|u_n(x)-E_xu_n(X_{t\wedge\tau_{D^m_{a,\varepsilon}}})|&\le c_2
a^\gamma\lim_{t\rightarrow 0^+}\sup_{n\ge
1}E_x\int_0^t\,dA_r^{\mu_n}\\& \le c_2a^\gamma\lim_{t\rightarrow
0^+}\sup_{n\ge 1}(\delta+E_xA_t^\mu) =c_2a^\gamma\delta.
\end{align*}
Since $\delta>0$ was arbitrary, we get
\begin{equation}
\label{eq1.28} \lim_{t\rightarrow 0^+}\sup_{n\ge1}
|u_n(x)-E_xu_n(X_{t\wedge\tau_{D^m_{a,\varepsilon}}})|=0.
\end{equation}
By the definition of the set $(E\setminus A_{a,\varepsilon}^m)^r$,
\begin{equation}
\label{eq1.29} P_x(\tau_{D^m_{a,\varepsilon}}>0)=1,\quad x\in
D^m_{a,\varepsilon}.
\end{equation}
By the Tanaka-Meyer formula and (\ref{eq1.2}), for every stopping
time $\tau$ we have
\begin{equation}
\label{eq4.ub} u_n^{\gamma+1}(X_\tau)\le
c_2(1+\gamma)E_x\Big(\int_0^\zeta\,dA_t^{\mu_n}|\FF_\tau\Big).
\end{equation}
It is clear that the family $\{A^{\mu_n}_\zeta\}$ is uniformly
integrable under $P_x$ for q.e. $x\in E$. Therefore the family
$\{E_x(\int_0^\zeta\,dA_t^{\mu_n}|\FF_\tau),n\ge 1\}$ is uniformly
integrable under $P_x$, and hence for fixed $\tau\in\mathcal{T}$
the family  $\{u_n(X_\tau),n\ge 1\}$ is uniformly integrable under
$P_x$ for q.e. $x\in E$.  From this and (\ref{eq1.29}) it follows
that for every $x\in D^m_{a,\varepsilon}$,
\[
\lim_{t\rightarrow 0^+}\sup_{n\ge 1}
|E_xu_n(X_{t\wedge\tau_{D^m_{a,\varepsilon}}})-
E_xu_n(X_t)\mathbf{1}_{\{t<\tau_{D^m_{a,\varepsilon}}\}}|\le
\lim_{t\rightarrow 0^+}\sup_{n\ge 1}
\int_{\{t\ge\tau_{D^m_{a,\varepsilon}}\}}
|u_n(X_{\tau_{D^m_{a,\varepsilon}}})|
=0.
\]
As a result,
\begin{equation}
\label{eq4.cp} \lim_{t\rightarrow 0}\sup_{n\ge 1}
|u_n(x)-p^{D^m_{a,\varepsilon}}_tu_n(x)|=0 \quad \mbox{q.e.  on
}D^m_{a,\varepsilon}.
\end{equation}
By (\ref{eq4.ub}),
\begin{equation}
\label{eq4.ub1} u^{\gamma+1}_n(x)\le c_2(1+\gamma)v_n(x)
\end{equation}
for q.e. $x\in E$. Since $\{v_n\}$  converges quasi-uniformly,
there exists a nest, and we may assume that it is $\{F_n\}$, such
that $\{u_n\}$ is uniformly bounded on $F_k$ for every $k\ge 1$.
Therefore by \cite[Theorem 2.2, Proposition 2.4]{Kl:MA1},
$\{u_n\}$ has a subsequence (still denoted by $(n)$) such that
$\{u_n\}$ is convergent and its limit is finite for q.e. $x\in
D^m_{a,\varepsilon}$. Let $a_n\nearrow\infty$ and let
$A_n=A_{a_n,(2 a_n)^{-1}}^n$\,, $D_n=D_{a_n,(2 a_n)^{-1}}^n$\,. By
$F$ let us denote the fine support of $\mu$. Since $\mu$ is
nontrivial, Cap$(F)>0$. Therefore by (\ref{eq5.strp}) there exist $n_0$  such
that Cap$(\eta^{n_0}_{n_0,\frac{1}{2a_{n_0}}},F)>0$. Since
$(\EE,D[\EE])$ satisfies ($\EE.7$),   we have
\[
\psi^n_{(2a_{n})^{-1}}>0,\quad n\ge n_0\quad \mbox{q.e}.
\]
Therefore, by (\ref{eq4.ab}) and Lemma \ref{lm3.cap11},
\[
\mbox{Cap}(\psi^n_{n,(2a_n)^{-1}}<a_n^{-1})\le
\mbox{Cap}(\psi^{n_0}_{n,(2 a_{n_0})^{-1}}<a^{-1}_n)\rightarrow 0,
\quad n\rightarrow \infty.
\]
Since $\{A_n\}$
is a nest, it follows that
\begin{equation}
\label{eq1.30} \lim_{n\rightarrow\infty}\mbox{Cap}(E\setminus
A_n)=\lim_{n\rightarrow\infty}E_\pi\int_{D_{E\setminus
A_n}}^\infty e^{-t}\varphi(X_t)\,dt=0,
\end{equation}
the first equality being a consequence of \cite[Theorem
IV.5.28]{MR}. By \cite[Theorem IV.5.28]{MR} again and \cite[Lemma
V.2.19]{MR},
\begin{equation}
\label{eq1.31} \mbox{Cap}(E\setminus D_n)=E_\pi\int_{D_{E\setminus
D_n}}^\infty e^{-t}\varphi(X_t)\,dt
\end{equation}
for $n\ge1$. Without loss of generality we may assume that the
sequence $\{A_n\}$ is increasing, and consequently that  $\{D_n\}$
is increasing, for  otherwise we can replace $\{A_n\}$ by
$\{\tilde{A}_n\}$, where $\tilde{A}_n=\bigcup_{k=1}^{n} A_k$, and
consider  $\tilde{D}_n=\tilde{A}_n\setminus(E\setminus
\tilde{A}_n)^r$ in place of $D_n$. Therefore by (\ref{eq1.30}),
$P_\pi(\lim_{n\rightarrow\infty}D_{E\setminus A_n}<\zeta)=0$.
Since for every $B\in\BB(E)$,
\begin{equation}
\label{eq1.32} \sigma_B=D_B\quad\mbox{on }\{D_B>0\},
\end{equation}
we deduce from (\ref{eq1.30}) that
\begin{equation}
\label{eq4.abc}
P_\pi(\lim_{n\rightarrow\infty}\tau_{A_n}<\zeta)=0.
\end{equation}
By \cite[Proposition 10.6]{Sharpe},
\[
\tau_{A_n}=\tau_{D_n},\quad P_\pi\mbox{-a.s.}
\]
Hence $P_\pi(\lim_{n\rightarrow\infty}\tau_{D_n}<\zeta)=0$ and by
(\ref{eq1.32}), $P_\pi(\lim_{n\rightarrow\infty}D_{E\setminus
D_n}<\zeta)=0$. This and  (\ref{eq1.31}) show that
\begin{equation}
\label{eq1.33} \lim_{n\rightarrow\infty}\mbox{Cap}(E\setminus
D_n)=0.
\end{equation}
We have proved that for every $m\ge1$ there exists a subsequence
(still denoted by $(n)$) such that $\{u_n\}$ converges q.e. and its
limit is finite q.e. on $D_m$. Therefore by (\ref{eq1.33}) one can
find a further subsequence (still denoted by $(n)$) such that
$\{u_n\}$ converges q.e. and its limit is finite q.e. on $E$. Let
$w=\lim_{n\rightarrow\infty}u_n$ q.e. Since $\{u_n\}$ is q.e.
convergent,
\[
\sup_{n\ge1}E_xA_\zeta^{\nu_n}=\sup_{n\ge1}u_n(x)<\infty.
\]
Therefore by \cite[Section 4]{Jakubowski2} the sequence
$\{u_n(X)\}$ is uniformly $S$-tight under $P_x$ for q.e. $x\in E$.
It is also clear that for every $t\ge 0$, $u_n(X_t)\rightarrow
w(X_t)$ in probability $P_x$ for q.e. $x\in E$. Therefore by
\cite[Theorem 1]{Jakubowski1}, the  definition of the sets
$\{A_{a,\varepsilon}^m\}$ and the Lebesgue dominated  convergence
theorem,
\begin{equation}
\label{eq1.34}
E_x\int_0^{\tau_{A^m_{a,\varepsilon}}}g(u_n)(X_t)\,dA_t^{\mu_n}
\rightarrow
E_x\int_0^{\tau_{A^m_{a,\varepsilon}}}g(w)(X_t)\,dA_t^{\mu}
\end{equation}
as $n\rightarrow\infty$ for  q.e. $x\in E$. Moreover, since
$u_n\rightarrow w$ q.e.,
\begin{equation}
\label{eq1.35} u_n(X_{\tau_{A^m_{a,\varepsilon}}}) \rightarrow
w(X_{\tau_{A^m_{a,\varepsilon}}}),\quad P_x\mbox{-a.s.}
\end{equation}
for q.e. $x\in E$. By the definition of a solution of
(\ref{eq1.27}),
\[
u_n(X_t)=E_xu_n(X_{\tau_{A^m_{a,\varepsilon}}})
+E_x\int_0^{\tau_{A^m_{a,\varepsilon}}}g(u_n)(X_t)\,dA_t^{\mu_n}
\]
for q.e. $x\in E$. By the above,  (\ref{eq1.34}) and
(\ref{eq1.35}),
\begin{equation}
\label{eq1.36} w(x)=E_xw(X_{\tau_{A^m_{a,\varepsilon}}})
+E_x\int_0^{\tau_{A^m_{a,\varepsilon}}}g(w)(X_t)\,dA_t^{\mu}
\end{equation}
for q.e. $x\in E$. By (\ref{eq4.ub1}), $w$ is a potential.
Therefore replacing $A^m_{a,\varepsilon}$ in (\ref{eq1.36}) by
$A_n$, letting $n\rightarrow \infty$ and using (\ref{eq4.abc}) we
obtain
\[
w(x)=E_x\int_0^\zeta g(w)(X_t)\,dA_t^\mu
\]
for q.e. $x\in E$. By uniqueness, $w=u$.
\end{dow}

\begin{stw}
\label{stw1.8} Let $\mu_n,\mu\in\MM_{0,b}$. If
$\|\mu_n-\mu\|_{TV}\rightarrow 0$ then $\mu_n\uAF\mu$.
\end{stw}
\begin{dow}
Let $u_n(x)=E_xA_\zeta^{\mu_n}$, $u(x)=E_xA_\zeta^{\mu}$.  By
\cite[Lemma 6.1]{BDHPS}, for every $q\in(0,1)$,
\[
E_x\sup_{t\ge 0}|u_n(X_t)-u(X_t)|^q
\le\frac{1}{1-q}\Big(E_x\int_0^\zeta\,dA_t^{|\mu_n-\mu|}\Big)^q
\]
for q.e. $x\in E$, where  $|\mu^n-\mu|$ stands for the total
variation of the measure $\mu_n-\mu$. Let $\beta\in S_{00}^{(0)}$
be such that $\beta(E)=1$. Then from the above inequality we
conclude that for every $q\in(0,1)$,
\begin{align*}
E_\beta\sup_{t\ge 0}|u_n(X_t)-u(X_t)|^q&
\le\frac{1}{1-q}\Big(E_\beta\int_0^\zeta\,dA_t^{|\mu_n-\mu|}\Big)^q
\le\frac{1}{1-q}\|R_\beta\|^q_\infty\|\mu-\mu_n\|^q_{TV}.
\end{align*}
By Remark \ref{uw1.1} there exists a subsequence (still denoted by
$(n)$) such that $u_n\rightarrow u$ quasi-uniformly. Therefore the
proposition follows from Proposition \ref{stw1.3}.
\end{dow}
\medskip

The following proposition  answers the question raised in
\cite[Remark 3.6]{BoccardoOrsina}.

\begin{tw}
\label{stw1.9} Assume that $g$ satisfies \mbox{\rm(H)}. Let
$\{\mu_n\}\subset\mathcal{R}$ be nontrivial, $\{\nu_n\}\subset\MM_{0,b}^+$ and
let $u_n,v_n$ denote solutions of the problems
\[
-Au_n=g(u_n)\cdot\mu_n,\quad u_n>0,\qquad
-Av_n=g(v_n)\cdot(\nu_n+\mu_n),\quad v_n>0.
\]
If $u_n\rightarrow 0$ in the topology of $m$-a.e. convergence,
$\sup_{n\ge1}\|\nu_n\|_{TV}<\infty$ and  $\nu_n\uAF\nu $ for some nontrivial
$\nu\in\MM_{0,b}^+$ then $v_n\rightarrow v$ in the topology of
$m$-a.e. convergence, where $v$ is a solution of
\begin{equation}
\label{eq1.44} -Av=g(v)\cdot\nu,\quad v>0.
\end{equation}
\end{tw}
\begin{dow}
By Proposition \ref{stw1.1}, $u_n\le v_n$, so by monotonicity of
$g$,
\[
E_x\int_0^\zeta g(v_n)(X_t)\,dA_t^{\mu_n} \le E_x\int_0^\zeta
g(u_n)(X_t)\,dA_t^{\mu_n}=u_n(x).
\]
By the assumptions of the proposition, up to a subsequence,
\begin{equation}
\label{eq1.45} E_x\int_0^\zeta
g(v_n)(X_t)\,dA_t^{\mu_n}\rightarrow0
\end{equation}
as $n\rightarrow\infty$ for $m$-a.e. $x\in E$. Let $w_n$ be a
solution of
\[
-Aw_n=g(w_n)\cdot\nu_n,\quad w_n>0.
\]
By the Meyer-Tanaka formula,
\begin{align*}
|w_n(x)-v_n(x)|&\le E_x\int_0^{\zeta}
\mbox{sgn}(w_n-v_n)(X_t)(g(w_n)(X_t)-g(v_n)(X_t))\,dA_t^{\nu_n}\\
&\quad-E_x\int_0^\zeta\mbox{sgn}(w_n-v_n)g(v_n)(X_t)\,dA_t^{\mu_n}\\
&\le E_x\int_0^\zeta g(v_n)(X_t)\,dA_t^{\mu_n}.
\end{align*}
By the above estimate and (\ref{eq1.45}), up to a subsequence we
have $|w_n-v_n|\rightarrow 0$, $m$-a.e. Since
$\sup_{n\ge1}\|\nu_n\|_{TV}<\infty$, applying Theorem \ref{stw1.6}
shows that, up to a subsequence, $w_n\rightarrow v$, $m$-a.e.,
which completes the proof.
\end{dow}

\begin{stw}
\label{stw1.10} Let $g$, $\{\mu_n\}$ satisfy the assumptions of
Theorem \ref{stw1.9} and $u_n, v_n$ be as in Theorem \ref{stw1.9}, with
$\{\nu_n\}\subset\MM_{0,b}^+$ such that $\|\nu_n-\nu\|_{TV}\rightarrow
0$ for some nontrivial $\nu\in\MM_{0,b}^+$. Then  $v_n\rightarrow v$ in the
topology of  $m$-a.e. convergence.
\end{stw}
\begin{dow}
Follows from Proposition \ref{stw1.8} and Theorem \ref{stw1.9}.
\end{dow}
\medskip


We close this section with some results showing that
$\{\mu_n\}\subset\MM_{0,b}$ is ``locally equidiffuse" if it
converges in the $uAF$ sense. These results will not be needed
later on in our  study of stability of solutions of
(\ref{eq1.17}). However, we find them interesting and we think
that they shed a new light on the nature of the convergence in the
$uAF$ sense.

Let us recall that a family $\{\mu_t,t\in T\}\subset\MM_{0,b}$ is
called equidiffuse if for every $\varepsilon>0$ there exists
$\delta>0$ such that for every $A\in\BB(E)$, if Cap$(A)<\delta$
then $|\mu_t|(A)<\varepsilon$ for every $t\in T$.

\begin{stw}
\label{stw1.12} Assume that $(\EE,D[\EE])$ satisfies $(\EE.7)$.
Let $\mu,\mu_n\in\MM_{0,b}$ be such that
$\sup_{n\ge1}\|\mu_n\|_{TV}<\infty$ and $\mu_n\uAF\mu$. Then there
exists a bounded excessive function $\eta\in D_e[\EE]$ such that
$\eta>0$ and the family $\{\eta\cdot\mu_n\}$ is equidiffuse.
\end{stw}
\begin{dow}
By Proposition \ref{stw1.11}, if $u_n, u$ are defined by
(\ref{eq4.iff}) then, up to a subsequence, $u_n\rightarrow u$
quasi-uniformly. It follows that there exists a nest $\{F_k\}$
such that for every $k\ge1$,
\begin{equation}
\label{eq1.47} \sup_{n\ge 1}\sup_{x\in F_k}u_n(x)<\infty.
\end{equation}
Since $m$ is a smooth measure, there exists a nest
$\{\tilde{F}_n\}$ such that
$\|R\mathbf{1}_{\tilde{F}_n}\|_\infty<\infty$ and $m(\tilde{F}_n)<\infty$ for $n\ge 1$.
Therefore there exists a closed set $F$ such that $\eta:=
R\mathbf{1}_{F}>0$, $\|R\mathbf{1}_F\|_\infty<\infty$,
$m(F)<\infty$ and $F\subset F_k$ for some $k\ge 1$. It is clear
that $\eta$ is  excessive and $\eta\in D_e[\EE]$. Let
$\beta:=\mathbf{1}_{F}\cdot m$. Then for every $B\in\BB(E)$,
\begin{equation}
\label{eq1.48} \int_B\eta\,d\mu_n=\int_BR\mathbf{1}_{F}\,d\mu_n
=E_\beta\int_0^\zeta\mathbf{1}_B(X_t)\,dA_t^{\mu_n} \le E_\beta
\int_{D_B}^\zeta dA_t^{\mu_n}.
\end{equation}
The family $\{A_\zeta^{\mu_n}\}$ is uniformly integrable under the
measure $P_\beta$. To see this, let us first observe that by
(\ref{eq1.47}),
\[
\sup_{n\ge 1}\int_E|u_n(x)|^2\,\beta(dx)
=\sup_{n\ge1}\int_{F}|u_n(x)|^2\,m(dx)<\infty.
\]
Since $u_n\rightarrow u$, $m$-a.e., it follows that
\[
E_\beta A_\zeta^{\mu_n}=\int_E u_n(x)\,\beta(dx) \rightarrow\int_E
u(x)\,\beta(dx)=E_\beta A_\zeta^\mu.
\]
On the other hand, since $\mu_n\uAF\mu$,
$A_\zeta^{\mu_n}\rightarrow A_\zeta^\mu$ in measure $P_\beta$,
which proves that $\{A_\zeta^{\mu_n}\}$ is uniformly integrable
under $P_\beta$. The uniform integrability implies that
\begin{equation}
\label{eq4.ll} \lim_{n\rightarrow\infty}E_\beta\sup_{t\ge0}
|A^{\mu_n}_t-A^\mu_t|=0.
\end{equation}
Suppose that, contrary to our claim, the family
$\{\eta\cdot\mu_n\}$ is not equidiffuse. Then there exist
$\varepsilon>0$ and a sequence $\{B_k\}$ of Borel subsets of $E$
such that Cap$(B_k)\rightarrow 0$ and $\sup_{n\ge 1}\int_{B_k}\eta
d\mu_n\ge\varepsilon$, $k\ge 1$. Then by Theorem IV.5.28 and Lemma
2.19 in \cite{MR},
\[
P_\beta(\lim_{k\rightarrow\infty}D_{B_k}\wedge\zeta=\zeta)=1.
\]
From this, (\ref{eq1.48}) and (\ref{eq4.ll}) it follows that
$\sup_{n\ge 1}\int_{B_k}\eta\,d\mu_n\rightarrow 0$ as
$k\rightarrow\infty$. This leads to the contradiction that
$\{\eta\cdot\mu_n\}$ is not equidiffuse.
\end{dow}

\begin{wn}
Let $\{\mu_n\}$ be as in Proposition \ref{stw1.12}. If
$(\EE,D[\EE])$ is strongly Feller then for every compact $K\subset
E$ the family $\{\mathbf{1}_K\cdot\mu_n\}$ is equidiffuse.
\end{wn}
\begin{dow}
Follows from the fact that every excessive function with respect
to a strongly Feller Dirichlet form is lower semi-continuous.
\end{dow}

\nsubsection{Stability: Approximation of measures by
mollification} \label{sec7}

In this section we assume that $\mu$ is a nontrivial Borel measure
on a subset $E$ of $\BR^d$. By putting $\mu(\BR^d\setminus E)=0$
we may and will assume that $\mu$ is a Borel measure on $\BR^d$.
We study stability of solutions $u_n$ of (\ref{eq1.17}) in case
\begin{equation}
\label{eq7.1} \mu_n=j_{1/n}\ast\mu,\quad n\ge1,
\end{equation}
where $j_\varepsilon(x)=\varepsilon^{-d}j(\varepsilon^{-1}x)$ for
$x\in\BR^d$, $\varepsilon>0$ and
\[
j(x)=\left\{
\begin{array}{l}c\exp(\frac{1}{|x|^2-1}),\,\quad |x|<1, \\
0,\qquad\qquad\qquad |x|\ge 1
\end{array}
\right.
\]
with $c>0$  chosen so that $\int_{\BR^d}j(x)\,dx=1$. By
(\ref{eq5.2}), $u_n$ is a solution of the equation
\begin{equation}
\label{eq7.2} -Au_n=g(u_n)\cdot(j_{1/n}*\mu_d+j_{1/n}*\mu_c).
\end{equation}
We shall show that for some class of operators Theorem
\ref{stw1.9} is applicable to (\ref{eq7.2}). To this end, we first
consider the case $\mu_d=0$ in Theorem \ref{tw1.6} below, and then
we show that $j_{1/n}\ast\mu_d\uAF\mu_d$.

In the proof of the following theorem a key role is played by
Theorem \ref{tw1.5}.

\begin{tw}
\label{tw1.6} Assume that $(\EE,D[\EE])$ is a form on $E\subset
\BR^{d}$ satisfying $(\EE.7)$ and Meyer's hypothesis
{\mbox{\rm(L)}}. Let $\mu\in\MM^+_b$ be a nontrivial measure such
that $\mu\bot\mbox{\rm Cap}$. Let $u_n$ denote a solution of the
problem
\[
-Au_n=g(u_n)\cdot\mu_n,\quad u_n>0
\]
with $\mu_n$ defined by \mbox{\rm(\ref{eq7.1})}. Then
$u_n\rightarrow 0$ in the topology of $m$-a.e. convergence as
$n\rightarrow\infty$.
\end{tw}
\begin{dow}
Let $B\in \mathcal{B}(E)$ be such that Cap$(B)=0$ and
$\mu(E\setminus B)=0$.  Since $\mu$ is finite, there exists an
increasing sequence $\{F_k\}$ of closed subsets of $E$ such that
$\mu(B\setminus\bigcup_{k=1}^\infty F_k)=0$. Let
$\mu^k=\mathbf{1}_{F_k}\cdot\mu$, $\mu^k_n=j_{1/n}\ast\mu^k$. Then
$\mu=\lim_{k\rightarrow\infty}\mu^k$ and
$\mu_n=\lim_{k\rightarrow\infty}\mu^k_n$ in the total variation
norm. Without loss of generality we may assume that
$\|\mu^k-\mu\|_{TV}\le k^{-1}$ for $k\ge 1$. Let
$\nu^k_n=g(u^k_n)\cdot\mu^k_n$, $k,n\ge1$, where $u^k_n$ is a
solution of
\[
-Au^k_n=g(u^k_n)\cdot\mu^k_n, \quad u^k_n>0.
\]
By Theorem \ref{tw1.5}, for every sequence $(n')$ there exists a
subsequence (still denoted by $(n')$) and a smooth measure $\nu^k$
such that $\nu^k_{n'}\rightarrow\nu^k$ vaguely and
$u^k_{n'}\rightarrow u^k$, $m$-a.e. as $n'\rightarrow\infty$,
where $-Au^k=\nu^k$. For a closed set $F\subset E$ and $n\ge1$
write $B(F,n)=\{x\in E:\mbox{dist}(x,F)\le1/n\}$. By the
properties of the vague convergence, for every $n\ge1$ we have
\[
0=\liminf_{n'\rightarrow\infty}\nu^k_{n'}(E\setminus B(F_k,n))\ge
\nu^k(E\setminus B(F_k,n)).
\]
Since this holds for every $n\ge1$, $\nu^k(E\setminus F_k)=0$. Hence
$\nu^k\equiv 0$, because Cap$(F_k)=0$ and $\nu^k$ is a smooth
measure. As a consequence, $u^k=0$. By Proposition \ref{stw1.1},
$u^k_n\le u_n$. By the Meyer-Tanaka formula, (H) and
(\ref{eq1.2}),
\begin{align*}
|u_n(x)-u^k_n(x)|^{\gamma+1}&\le(1+\gamma)E_x
\int_0^\zeta(u_n-u^k_n)^{\gamma}(X_t) (g(u_n)(X_t)\,dA_t^{\mu_n}
-g(u^k_n)(X_t)\,dA_t^{\mu^k_n})
\\&=(1+\gamma)E_x
\int_0^\zeta(u_n-u^k_n)^\gamma(X_t)
(g(u_n)-g(u^k_n))(X_t)\,dA_t^{\mu^k_n}\\
&\quad+(1+\gamma)E_x \int_0^\zeta g(u_n)(X_t)
(u_n-u^k_n)^{\gamma}(X_t)\,(dA_t^{\mu_n}
-dA_t^{\mu^k_n})\\
&\le(1+\gamma)E_x \int_0^\zeta g(u_n)(u_n)^\gamma(X_t)\,
(dA_t^{\mu_n}-dA_t^{\mu^k_n})\\
&\le(1+\gamma)c_2 E_x\int_0^\zeta\,dA_t^{|\mu_n-\mu^k_n|}.
\end{align*}
Let $\beta\in S_{00}^{(0)}$. From the above inequality we conclude
that
\begin{align*}
\int_E|u_n-u^k_n|^{1+\gamma}\,d\beta\le(1+\gamma)
c_2E_\beta\int_0^\zeta\,dA_t^{|\mu_n-\mu^k_n|} &\le (1+\gamma)c_2
\|R\beta\|_\infty\cdot\|\mu_n-\mu_n^k\|_{TV}\\
&\le (1+\gamma)c_2\|R\beta\|_\infty k^{-1}.
\end{align*}
Therefore for every  $\beta\in S^{(0)}_{00}$,
\[
\int_Eu_{n'}\,d\beta \le\int_E|u_{n'}-u^k_{n'}|\,d\beta
+\int_Eu^k_{n'}\,d\beta \le
c(\beta,\gamma,c_2)k^{-1}+\int_Eu^k_{n'}\,d\beta.
\]
Letting $n'\rightarrow\infty$ and then $k\rightarrow\infty$ in the
above inequality we see that $\int_Eu_{n'}\,d\beta\rightarrow 0$
for every finite $\beta\in S_{00}^{(0)}$, which implies that, up
to a subsequence, $u_{n'}\rightarrow 0$, $m$-a.e.
\end{dow}
\medskip

In the rest of the section we confine ourselves to the class of
forms defined below. Let $\psi:\BR^d\rightarrow\BR$ be defined as
\begin{equation} \label{eq4.ls} \psi(x)=\frac12
(Bx,x)+\int_{\BR^d}(1-\cos(x,y))J(dy),
\end{equation}
where $B$ is a $d$-dimensional nonnegative definite symmetric
matrix and $J$ is a symmetric Borel measure on
$\BR^d\setminus\{0\}$ satisfying
\[
\int_{\BR^d\setminus\{0\}}\frac{|x|^2}{1+|x|^2}J(dx)<\infty.
\]
Consider the form $(\BB,D[\BB])$ on $L^2(\BR^d;dx)$ defined as
\begin{equation}
\label{eq4.form}
\left\{
\begin{array}{l}\BB(u,v)=\int_{\BR^d}\hat{u}(x)\bar{\hat{v}}(x)\psi(x)\,dx,
\smallskip\\
D[\BB]=\{u\in L^2(E;dx);\, \int_{\BR^d}|\hat{u}(x)|^2\psi(x)\,dx<\infty\},
\end{array}
\right.
\end{equation}
where $\hat{u}$ stands for the Fourier transform of $u$. It is
well known (see \cite[Example 1.4.1]{Fukushima}) that
$(\BB,D[\BB])$ is a symmetric regular Dirichlet form on
$L^2(\BR^d;dx)$.

An important example of  $\psi$ of the form (\ref{eq4.ls}) is
$\psi(x)=\phi(|x|^2)$, $x\in\BR^d$, where
$\phi:(0,\infty)\rightarrow [0,\infty)$ is a Bernstein function,
i.e. smooth function such that $(-1)^n D^n\phi\le 0$ for $n\ge1$.
In this case the operator $A$ associated with the form
$(\BB,D[\BB])$ is equal to $\phi(\Delta)$. For instance,
$A=\Delta^{\alpha/2}$ for $\phi(x)=x^{\alpha/2}$ with $\alpha\in
(0,2]$.

It is  well known (see \cite[Example 1.4.1]{Fukushima}) that if
$(\BB,D[\BB])$ satisfies Meyer's hypotheses (L) then the
$\alpha$-Green function $G_\alpha(\cdot,\cdot)$ has the property
that
\[
G_\alpha(x,y)=G_\alpha(x-y),\quad x,y\in\BR^d
\]
for some real function $G_\alpha$ defined on $\BR^d$.

For an arbitrary open set $D\subset \BR^d$ let $(\EE,D[\EE])$
denote the part of $(\BB,D[\BB])$ on $D$, i.e.
\begin{equation}
\label{eq4.dfp} D[\EE]=\{u\in D[\BB]:u= 0,\,\, m\mbox{-a.e  on }
\BR^d\setminus D\},\quad \EE(u,v)=\BB(u,v),\, u,v\in D[\EE].
\end{equation}
By \cite[Theorems 4.4.3, 4.4.4]{Fukushima}, $(\EE,D[\EE])$ is a
symmetric regular Dirichlet form on $L^2(D;dx)$. For instance, if
$(\BB,D[\BB])$ is defined by (\ref{eq4.form}) with
$\psi(x)=|x|^{\alpha}$ for some $\alpha\in(0,2]$ and $A_D$ is the
operator associated with $(\EE,D[\EE])$ then the solution of the
problem
\[
-A_Du=f
\]
with $f\in L^2(D;dx)$ may be interpreted  as  a solution of the
Dirichlet problem
\[
-\Delta^\alpha u=f\mbox{ in }D,\quad u(x)=0\mbox{ for }x\in
\BR^d\setminus D.\quad
\]


Let us recall that in the whole paper we assume that the forms
under consideration satisfy ($\EE.5$), ($\EE.6$). We have already
mentioned that the form defined  by (\ref{eq4.dfp}) satisfies
($\EE.5$), i.e. is regular. It is known, that it satisfies
($\EE.6$) if $\psi^{-1}$ is locally integrable on $\BR^d$ (see
\cite[Example 1.5.2]{Fukushima}).

For a given Hunt process $\mathbb{X}$ on $E$ and an open set
$D\subset E$  we denote by $\mathbb{X}^D$ the Hunt process on $D$
which is a part of the process $\mathbb{X}$ on $D$ (see
\cite[Appendix A.2]{Fukushima} for details).

\begin{stw}
\label{stw1.7} Assume that $(\EE,D[\EE])$  defined by
\mbox{\rm(\ref{eq4.dfp})} satisfies Meyer's hypothesis
\mbox{\rm{(L)}}.  Let  $\mu\in\MM_{0,b}$ and $\mu_n$ be defined by
\mbox{\rm(\ref{eq7.1})}. Then $\mu_n\uAF \mu$ as
$n\rightarrow\infty$.
\end{stw}
\begin{dow}
Let us first observe that the proof can be reduced to the case
$E=\BR^d$. This follows from the fact that if $\mathbb{X}$ is a
Hunt process associated with the form $(\EE,D[\EE])$ on $E=\BR^d$
and if $A^\mu$ is an additive
functional of $\mathbb{X}^D$ associated with a smooth measure
$\mu$ on $D$ then
\[
A^\mu_t=A^{\bar{\mu}}_{t\wedge \tau_D}, \quad t\ge 0,
\]
where $A^{\bar{\mu}}$ is the additive functional of $\mathbb{X}$
associated with the measure $\bar{\mu}$ being the extension of
$\mu$ to $\BR^d$ by putting zero on $\BR^d\setminus D$.

Let $u\in D_e[\EE]$ and let $u_\varepsilon=j_\varepsilon\ast u$.
Then
\begin{align}
\label{eq1.37} \nonumber\EE(u_\varepsilon,u_\varepsilon)
=\int_{\BR^d}\hat{u}_{\varepsilon}(x)\cdot\bar{\hat{u}}_{\varepsilon}(x)
\cdot\psi(x)\,dx&= \int_{\BR^d}|\hat{u}|^2(x)\cdot|\hat{
j_\varepsilon}|^2(x)\psi(x)\,dx\\&\le
\int_{\BR^d}|\hat{u}|^2(x)\psi(x)\,dx=\EE(u,u).
\end{align}
Observe that for every $\alpha\ge 0$,
\begin{equation}
\label{eq1.37a}
R_\alpha u(x)=\int_{\BR^d}G_\alpha(x-y)u(y)\,dy=
(G_\alpha\ast u)(x).
\end{equation}
Hence
\[
R_\alpha u_\varepsilon=G_\alpha\ast u_\varepsilon
=G_\alpha\ast(u\ast j_\varepsilon)= j_\varepsilon\ast(G_\alpha\ast
u)=j_\varepsilon\ast(R_\alpha u).
\]
In particular, for every $u\in D(A)$, $j_\varepsilon\ast u\in
D(A)$ and
\begin{equation}
\label{eq1.38}
-A(j_\varepsilon\ast u)=j_\varepsilon\ast(-Au).
\end{equation}
Assume that $u\in D(A)$ and write  $u_n=j_{1/n}\ast u$. Applying
(\ref{eq1.38}) gives
\begin{align*}
\|u_n-u_m\|_\EE=(-A(u_n-u_m),u_n-u_m)&
=(j_{1/n}\ast(-Au)-j_{1/m}\ast(-Au),u_n-u_m)\\
&\le2\|-Au\|_{L^2}\|u_n-u_m\|_{L^2}.
\end{align*}
Hence $u_n\rightarrow u$ in $\EE$. Now assume that $u\in D[\EE]$.
Then by (\ref{eq1.37}),
\begin{align*}
\|u_n-u\|_\EE&\le\|u_n-j_{1/n}\ast(\alpha R_\alpha u)\|_\EE
+\|j_{1/n}\ast(\alpha R_\alpha u)-\alpha R_\alpha u\|_\EE
+\|\alpha R_\alpha u-u\|_\EE\\
& \le 2\|\alpha R_\alpha u-u\|_\EE +\|j_{1/n}\ast(\alpha R_\alpha
u)-\alpha R_\alpha u\|_\EE.
\end{align*}
Letting $n\rightarrow\infty$ and then $\alpha\rightarrow\infty$ we
conclude from the above inequality that $j_{1/n}\ast u\rightarrow
u$ in $\EE$. Finally, assume that $u\in D_e[\EE]$. Then there
exists a sequence $\{u^k\}\subset D[\EE]$ such that
$\|u^k-u\|_\EE\rightarrow 0$. Using once again (\ref{eq1.37}) we
obtain
\begin{align*}
\|u_n-u\|_\EE&=\|u_n-j_{1/n}\ast u^k\|_\EE +\|j_{1/n}\ast
u^k-u^k\|_\EE +\|u-u^k\|_\EE\\
& \le 2\|u-u^k\|_\EE+\|j_{1/n}\ast u^k-u^k\|_\EE.
\end{align*}
Letting $n\rightarrow\infty$ and then $k\rightarrow\infty$ shows
that  $j_{1/n}\ast u\rightarrow u$ in $\EE$. Accordingly, for
every $u\in D_e[\EE]$,
\begin{equation}
\label{eq1.39} j_{1/n}\ast u\rightarrow u\quad\mbox{in }\EE.
\end{equation}
Let $\mu\in S_0^{(0)}$ and let $u_n,u$ be solutions of the
problems
\begin{equation}
\label{eq1.40a} -Au_n=j_{1/n}\ast\mu,\quad -Au=\mu.
\end{equation}
Then by (\ref{eq1.37a}),
\[
u=R\mu,\quad u_n=R(j_{1/n}\ast\mu)=j_{1/n}\ast(R\mu).
\]
Since $\mu\in S_0^{(0)}$, $R\mu\in D_e[\EE]$. Therefore by
(\ref{eq1.39}),
\begin{equation}
\label{eq1.40} \lim_{n\rightarrow\infty}\|u_n-u\|_\EE=0.
\end{equation}
From this and \cite[Lemma 5.1.1]{Fukushima} it follows that for
every $\alpha\in(0,1)$ and $\beta\in S_0^{(0)}$,
\begin{equation}
\label{eq1.41}
\lim_{n\rightarrow\infty}E_\beta\sup_{t\ge0}|u_n(X_t)-u(X_t)|^\alpha
\le c(\alpha,\beta)\lim_{n\rightarrow\infty}\|u_n-u\|_\EE=0.
\end{equation}
Now let  $u,u_n$ be solutions of (\ref{eq1.40a}) with
$\mu\in\MM_{0,b}$. Let $\{F_k\}$ be a nest such that
$\mathbf{1}_{F_k}\cdot\mu\in S_0^{(0)}$ for $k\ge1$, and let
$u^k_n, u^k$ be solutions of
\begin{equation}
\label{eq1.42} -Au_n^k=j_{1/n}\ast(\mathbf{1}_{F_k}\cdot\mu),\quad
-Au^k=\mathbf{1}_{F_k}\cdot\mu.
\end{equation}
From the probabilistic interpretation of equations
(\ref{eq1.40a}), (\ref{eq1.42}) and calculations  leading to
(\ref{eq1.20}) it follows that for every $\beta\in S_{00}^{(0)}$
such that $\beta(E)=1$ and every $q\in (0,1)$,
\begin{align}
\label{eq1.43} \nonumber &E_\beta\sup_{t\ge
0}|u_n^k(X_t)-u_n(X_t)|^q \le c(q,\beta)\|\mu^k-\mu\|_{TV},\,
\\&E_\beta\sup_{t\ge 0}|u(X_t)-u^k(X_t)|^q\le
c(q,\beta)\|\mu^k-\mu\|_{TV},\qquad
\end{align}
where $\mu^k=\mathbf{1}_{F_k}\cdot\mu$. For $u\in\BB(E)$ put
$|u|^{[q]}_{\sup}=E_\beta\sup_{t\ge 0}|u(X_t)|^q$. Then by
(\ref{eq1.41}) and (\ref{eq1.43}),
\begin{align*}
\lim_{n\rightarrow\infty}|u_n-u|^{[q]}_{\sup}&
\le\lim_{n\rightarrow\infty}(|u_n-u_n^k|^{[q]}_{\sup}
+|u_n^k-u^k|^{[q]}_{\sup} +|u^k-u|^{[q]}_{\sup})\\
&\le 2
c(q,\beta)\|\mu^k-\mu\|_{TV}
+\lim_{n\rightarrow\infty}|u_n^k-u^k|^{[q]}_{\sup}\\
&=2c(q,\beta)\|\mu^k-\mu\|_{TV}.
\end{align*}
Letting $k\rightarrow\infty$ shows that
$|u_n-u|^{[q]}_{\sup}\rightarrow 0$. The desired result now
follows from Proposition \ref{stw1.3}.
\end{dow}
\medskip

Combining Theorem \ref{stw1.9} with Theorems \ref{tw1.6} and
\ref{stw1.7} we get the following stability result.

\begin{tw}
\label{tw1.7} Let $(\EE,D[\EE])$ be the form defined by
{\mbox{\rm(\ref{eq4.dfp})}} such that  $(\EE.7)$  and Meyer's
hypothesis {\mbox{\rm(L)}} are satisfied. Assume that $g$
satisfies \mbox{\rm(H)}, $\mu\in\MM_b$, $\mu_d$ is nontrivial and
by $u_n,u$ denote solutions of the problems
\[
-Au_n=g(u_n)\cdot(j_{1/n}\ast\mu),\quad u_n>0,
\]
\begin{equation}
\label{eq7.15} -Au=g(u)\cdot\mu_d,\quad u>0.
\end{equation}
Then $u_n\rightarrow u$ in the topology of  $m$-a.e. convergence
as $n\rightarrow\infty$.
\end{tw}
\begin{dow}
Follows from Theorem \ref{stw1.9} applied to the sequences
$\{\mu_n=j_{1/n}\ast\mu_c\}$ and $\{\nu_n=j_{1/n}\ast\mu_d\}$. The
assumptions of Theorem \ref{stw1.9} for $\{\mu_n\}$ are satisfied
by Theorem \ref{tw1.6}, whereas the assumptions for $\{\nu_n\}$
are satisfied thanks to Proposition \ref{stw1.7}.
\end{dow}
\medskip

We see that in the limit equation (\ref{eq7.15}) the whole
concentrated  part of $\mu$ disappear.  This and the fact that
(\ref{eq7.15}) has a unique solution makes is legitimate to call
$u$ satisfying (\ref{eq7.15}) the solution of (\ref{eq1.1}). With
this definition in mind, Theorem \ref{tw1.7} may be viewed as an
existence theorem for equation (\ref{eq1.1}).

\vspace{3mm}

\noindent{\bf\large Acknowledgements}
\medskip\\
Research supported by National Science Centre Grant No.
2012/07/B/ST1/03508.

\end{document}